\documentclass[12pt]{iopart}

\usepackage{iopams}
\usepackage{graphicx}
\usepackage{float}

\usepackage{cite}

\usepackage{color}

\usepackage[unicode,bookmarks,bookmarksopen,bookmarksopenlevel=0,colorlinks,
pageanchor=true,breaklinks=true,linkcolor=red,citecolor=blue,urlcolor=red,
pdfencoding = auto, bookmarksdepth = 4]{hyperref}

\begin{document}

\title[Interactions between discontinuities and hodograph method]{Interactions between discontinuities for binary mixture separation
problem and hodograph method}

\author{M S Elaeva$^1$, E V Shiryaeva$^2$, M Yu Zhukov$^{3,\,4}$}

\address{$^1$ Faculty of Applied Mathematics and Computer Science, Financial University under the Government of the Russian Federation, 49 Leningradsky prospekt, 125993, Moscow, Russia}
\ead{mselaeva@fa.ru}

\address{$^2$ Institute of Mathematics, Mechanics and Computer Science, Southern Federal University, 105/42 Bolshaya Sadovaya Street, 344006, Rostov-on-Don, Russia}
\ead{shir@math.sdedu.ru}

\address{$^3$ Institute of Mathematics, Mechanics and Computer Science, Southern Federal University, 105/42 Bolshaya Sadovaya Street, 344006, Rostov-on-Don, Russia}
\ead{myuzhukov@gmail.com}

\address{$^4$ South Mathematical Institute, Vladikavkaz Center of RAS, 22 Markus Street, 362027, Vladikavkaz, Russia}

\vspace{10pt}

\begin{indented}
\item[]February 2016
\end{indented}

\begin{abstract}
The Cauchy problem for first-order PDE with the initial  data which have a piecewise discontinuities localized in different spatial points is completely solved. The  interactions between  discontinuities arising after breakup of initial discontinuities are studied with the help of the hodograph method. The solution is constructed in analytical implicit form. To recovery the explicit form of solution we propose   the transformation of the PDEs into some ODEs on the level lines (isochrones) of implicit solution. In particular, this method allows us to solve the Goursat problem with initial data on characteristics.
The paper describes a specific problem for zone electrophoresis (method of the mixture separation). However, the method proposed allows to solve any system of two first-order quasilinear  PDEs for which the second order linear PDE, arising after the hodograph transformation, has  the Riemann-Green function in explicit form.
\end{abstract}

\pacs{02.30.Ik, 02.30.Jr, 82.80.Bg, 87.15.Tt}

\ams{35Lxx, 35L67, 35L40, 35L45 35L50, 35L65}

%
%
%
%
%
%
%
%
%

		
%
%

\vspace{2pc}
\noindent{\it Keywords}: hyperbolic quasilinear equations, hodograph method, zone electrophoresis

\submitto{\jpa}

\maketitle

%



\section{Introduction}\label{zhshel:sec:1}

The paper is devoted to the investigation of the discontinuities interactions for the first-order quasilinear hyperbolic PDEs. We consider the Cauchy problem for systems of first-order hyperbolic quasilinear equations with discontinuous piecewise constant initial data and additional the Rankine--Hugoniot  conditions on the discontinuities 
\begin{eqnarray}\label{zhshel:eq:0.01}
  u^k_t + \varphi^k_x(u^1,u^2) = 0, \quad  u^k(x,0)=f^k(x), \quad D[u^k]=[\varphi^k], \quad k=1,2,
\end{eqnarray}
where $f^k(x)$ is the piecewise constant function, $D$ is the discontinuity velocity, $\left[ \,\cdot\, \right]$ is the notation for jump of function across discontinuity.


When the  discontinuities of the functions $f^k(x)$ are localized
only at one point $x_0$ the problem (\ref{zhshel:eq:0.01})  is well-known as the Riemann problem. In this case,
the Cauchy problem solution is a moving sequence of  shock waves
(strong discontinuities $S_{x_0}$) and rarefaction waves (weak discontinuities $W_{x_0}$).


We discuss a slightly more general problem in the case when $f^k(x)$ have discontinuities at two points $x_1$ and $x_2$.
Obviously, in this case, a moving shock wave and rarefaction wave arising at the points $x_1$ and $x_2$
are interacting with each other in the process of the solution evolution.
Four types of interactions are possible: (i) a shock wave $S_{x_1}$ with a shock wave $S_{x_2}$, (ii) a shock wave $S_{x_1}$ with a rarefaction wave $W_{x_2}$,
(iii) a rarefaction wave $W_{x_1}$ with a shock wave $S_{x_2}$, (iv) a rarefaction wave $W_{x_1}$ with a rarefaction wave $W_{x_2}$.


Studies have shown, that the case of interaction of two rarefaction waves is the most difficult to describe. Indeed,
the interaction of weak discontinuities there is the Goursat problem with the initial  data on the characteristics because the fronts of the rarefaction waves (weak discontinuities) move along the characteristics. In the case of a single discontinuity at the point $x_0$ the different ways of constructing the Goursat problem solution (and the Riemann problem) for two quasilinear equations are presented, for instance, in~\cite{RozhdestvenskiiYanenko}. These methods in general case are based on iterative method of the solution constructing. For one particular quasilinear equations, e.g., equations of isotachophoresis and chromatography, it is possible to obtain an analytical solution in the case when there are several discontinuities of the function $f(x)$. In the case of two
quasilinear equations, the problem of discontinuities interaction ($S_{x_1}, S_{x_2}$), ($S_{x_1}, W_{x_2}$), ($W_{x_1}, S_{x_2}$) for the equations of electrophoresis is completely solved in \cite{ElaevaIzvestiya}. To date, the problem of interaction type ($W_{x_1}, W_{x_2}$) was an open problem.


To solve problem of interaction between discontinuities we propose to use the hodograph method for system of two quasilinear hyperbolic equations written in the Riemann invariants. We recall that in this method the dependent and independent variables reverse roles ($(u^1,u^2) \rightleftharpoons (t,x)$). In this case, the solution is an implicit form $t=t(u^1,u^2)$, $x=x(u^1,u^2)$.
To determine, for example, the function $t(u^1,u^2)$ there is a system of first-order linear PDEs with variable coefficients (see, e.g., \cite{RozhdestvenskiiYanenko}). The solvability condition allows us to transform the system into a single second-order linear PDE with variable coefficients. For this equation the Riemann--Green function is constructed (see, e.g., \cite{Copson}) and then the Goursat problem solution is obtained (see e.g. \cite{Bizadze}).


In a recent paper \cite{SenashovYakhno} a variant of the hodograph method based on the presence of conservation laws for the Riemann invariants is proposed. This method allows  us to obtain effectively  an implicit solution in the form $t=t(u^1,u^2)$, $x=x(u^1,u^2)$ for the case when there is an explicit form of the Riemann--Green function, for example, for the shallow water equations,
chromatography equations (see e.g. \cite{RozhdestvenskiiYanenko,FerapontovTsarev_MatModel}), zone electrophoresis equations (see e.g. \cite{BabskiiZhukovYudovichRussian}), etc. Moreover, the solution of the Cauchy problem can be obtained for arbitrary initial data including discontinuities.


In the papers \cite{Zhuk_Shir_ArXiv_2014_Part2} the method \cite{SenashovYakhno} was only slightly modified. In addition, the  numerically-analytical method proposed allows to reconstruct the explicit solution $u^1=u^1(x,t)$, $u^2=u^2(x,t)$ with the help of its implicit form $t=t(u^1,u^2)$, $x=x(u^1,u^2)$. This allows us to transform the Cauchy problem for PDE to the Cauchy problem for ODEs on the level lines of an implicit solution. In fact, the calculations are required only for integration of the Cauchy problem for ODEs. In this paper the `numerical part' of the method is almost not required. Most of the relations describing the solution are obtained in an explicit analytical form.


We note that the problem of the discontinuities interaction has important practical application, in particular, for the method of electrophoresis. Electrophoresis is a separation method of the multicomponent mixture into individual components under action of an electric field. A diffusionless variant of the electrophoresis mathematical models  is the conservation laws for hyperbolic equations. The discontinuous initial  data  correspond to the boundaries of the spatial regions occupied by the individual components. The study of the discontinuities behavior  and their interactions allows us to understand the complex scenario of components motion for  the  separation  process.


The general formulation of the zone electrophoresis problem is given in~\cite{BabskiiZhukovYudovichRussian}.
In the case of the two-component mixture the interactions ($S_{x_1}, S_{x_2}$), ($S_{x_1}, W_{x_2}$), ($W_{x_1}, S_{x_2}$) are studied in
\cite{ElaevaIzvestiya,ElaevaMM,Elaeva_Diss,Elaeva_ZhVM}.
The results presented below indicate the analytical method of the solving for the problem of the rarefaction waves interaction.
This gives us a complete analytical solution of the problem (together with the results of~\cite{ElaevaIzvestiya,ElaevaMM,Elaeva_Diss,Elaeva_ZhVM}).


Note that the  method proposed is weak linked  with the electrophoresis problem.
This method can be used for any system of two hyperbolic quasilinear equations.
However, the application of the method requires explicit analytic form of the Riemann--Green function for the linear equation arising from the use of the hodograph method.


The paper is organized as follows. In section~\ref{zhshel:sec:2} the basic electrophoresis equations, written in the original variables (concentrations of components) and the Riemann invariants, are presented. In section~\ref{zhshel:sec:3} we formulate the Riemann problem for the case of initial data with discontinuities at different  points in space. The Cauchy problem with initial data in the moment of interaction of weak discontinuities is presented. In section~\ref{zhshel:sec:4} the hodograph method is described. In section~\ref{zhshel:sec:5} the method of recovering the explicit form of the solution from its implicit form is described. In sections~\ref{zhshel:sec:6}--\ref{zhshel:sec:8} different variants of the discontinuities interactions are considered. In \ref{zhshel:sec:9} the brief description of numerical-analytical method for the solving of the Cauchy problem is presented.




\section{Zone Electrophoresis Equations}\label{zhshel:sec:2}

To describe zone electrophoresis of binary mixture in diffusionless approximation we use two first-order quasilinear hyperbolic equations \cite{BabskiiZhukovYudovichRussian,ElaevaIzvestiya,ElaevaMM,Elaeva_Diss,Elaeva_ZhVM}
\begin{eqnarray}\label{zhshel:eq:2.01}
  \frac{\partial u^k}{\partial t} + \frac{\partial }{\partial x} \left( \frac{\mu^k u^k}{1+s} \right) = 0,
  \quad k=1,2,  \quad s = u^1+u^2 > -1,
\end{eqnarray}
where $u^k$ are the effective concentrations, $\mu^k$ are the component mobilities,
($1+s$) is the conductivity of the mixture, $E=(1+s)^{-1}$ is the intensity of the electric field.

System (\ref{zhshel:eq:2.01}) can be written in the Riemann invariants
\begin{eqnarray}\label{zhshel:eq:2.03}
  R^k_t + \lambda^k(R^1,R^2)R^k_x = 0, \quad
  \quad k=1,2,
\end{eqnarray}
\begin{eqnarray}\label{zhshel:eq:2.04}
  \lambda^1(R^1,R^2)=R^1 R^1 R^2, \quad   \lambda^2(R^1,R^2)=R^2 R^1 R^2.
\end{eqnarray}

The concentrations $u^1$, $u^2$ are connected to the Riemann invariants $R^1$, $R^2$ as
\begin{eqnarray}\label{zhshel:eq:2.05}
u^1=\frac{\mu^2 (R^1-\mu^1)(R^2-\mu^1)}{R^1 R^2(\mu^1-\mu^2)},\quad
u^2=\frac{\mu^1 (R^1-\mu^2)(R^2-\mu^2)}{R^1 R^2(\mu^2-\mu^1)}.
\end{eqnarray}

The roots of the  polynomial
\begin{eqnarray*}
P(R^1,R^2) \equiv (1+u^1+u^2)(R)^2-(\mu^1+\mu^2+u^1\mu^2+u^2\mu^1)R+\mu^1\mu^2
\end{eqnarray*}
allows us to get inverse relations $R^1=R^1(u^1,u^2)$, $R^2=R^2(u^1,u^2)$.



\section{The analogue of the Riemann problem}\label{zhshel:sec:3}

We study the equations (\ref{zhshel:eq:2.03}), (\ref{zhshel:eq:2.04}) with
the discontinuous initial data
\begin{eqnarray}\label{zhshel:eq:3.01}
R^1\bigr|_{t=0}=R^1_0(x), \quad R^2\bigr|_{t=0}=R^2_0(x),
\end{eqnarray}
where
\begin{eqnarray}\label{zhshel:eq:3.02}
R^1_0(x)=
\left\{
  \begin{array}{lll}
    \mu^1, & \quad x < x^1, \\
    q^1, & \quad x^1 < x < x^2, \\
    \mu^1, & \quad  x^2 < x,
  \end{array}
\right.
\quad
R^2_0(x)=
\left\{
  \begin{array}{lll}
    \mu^2, & \quad x < x^1, \\
    q^2, & \quad x^1 < x < x^2, \\
    \mu^2, & \quad  x^2 < x.
  \end{array}
\right.
\end{eqnarray}
Here, $q^1$, $q^2$ are given constants such that
\begin{eqnarray}\label{zhshel:eq:3.03}
0< q^1 < \mu^1 < \mu^2 < q^2.
\end{eqnarray}
Note that for the classical Riemann problem the initial data usually have the discontinuities at only one point (not $x_1$ and $x_2$).
We also note that for original equations (\ref{zhshel:eq:2.01}) the initial data (\ref{zhshel:eq:3.01}), (\ref{zhshel:eq:3.02}) correspond to zero concentrations outside the interval $(x^1,x^2)$. Inside the interval $(x^1,x^2)$ the initial concentrations are constant $u^1=u^1_0$, $u^2=u^2_0$. The values $u^1_0$, $u^2_0$ are determined with the help $q^1$, $q^2$, and (\ref{zhshel:eq:2.05}). Inequality (\ref{zhshel:eq:3.03}) means that the effective concentrations  $u^1_0 > 0$ and  $u^2_0<0$.

To solve the Riemann problem (\ref{zhshel:eq:2.03}), (\ref{zhshel:eq:2.04}), (\ref{zhshel:eq:3.01})--(\ref{zhshel:eq:3.03}) we add
the Rankine--Hugoniot conditions which correspond to conservation laws (\ref{zhshel:eq:2.01}).
Taking into account (\ref{zhshel:eq:2.05}) we rewrite these conditions for the Riemann invariants in the following form
\begin{eqnarray}\label{zhshel:eq:3.04}
D\left[\frac{(\mu^k-R^1)(\mu^k-R^2)}{\mu^k R^1 R^2}\right]=\left[(\mu^k-R^1)(\mu^k-R^2)\right],
\quad
k=1,2,
\end{eqnarray}
where $D$ is the discontinuity velocity, $\left[ \,\cdot\, \right]$ is the notation for jump of function across discontinuity.

In the case when the solution has strong discontinuities such solution should additionally satisfy to the Lax conditions.
These conditions are well known and we omit them (see e.g. \cite{Lax,RozhdestvenskiiYanenko,Liu} and also \cite{BabskiiZhukovYudovichRussian,ElaevaIzvestiya,ElaevaMM,Elaeva_Diss,Elaeva_ZhVM}).

For completeness, we present the results obtained in \cite{ElaevaIzvestiya,ElaevaMM,Elaeva_Diss,Elaeva_ZhVM}).
The equations (\ref{zhshel:eq:2.01})  (and (\ref{zhshel:eq:2.03}), (\ref{zhshel:eq:2.04})) have the hyperbolic type. In particular, it means that the breakup of the initial discontinuities at the points $x_1$ and $x_2$ can be analyzed independently.
In the initial moment $t=+0$ at the point $x_1$ and also at point $x_2$ arise the rarefaction waves (self-similar solutions, weak discontinuities) and shock waves (discontinuous solutions, strong discontinuities). After a time, strong and weak discontinuities of the solution arising at point $x_1$ interact with strong and weak discontinuities of the solution arising at point $x_2$.

We illustrate the evolution of solution using $(x,t)$-plane and zones on this plane. It is convenient to introduce the notation
\begin{eqnarray*} 
\mathbb{Z}_k=\{R^1_k(x,t),R^2_k(x,t);a_k(t),b_k(t)\},
\end{eqnarray*}
where $\mathbb{Z}_k$ is the zone on $(x,t)$-plane, $a_k(t)$, $b_k(t)$ are the zone boundaries at the moment $t$, $R^1_k(x,t)$ and $R^2_k(x,t)$ are the Riemann invariants in this zone.

\begin{figure}[H]
  \centering
\includegraphics[scale=0.95]{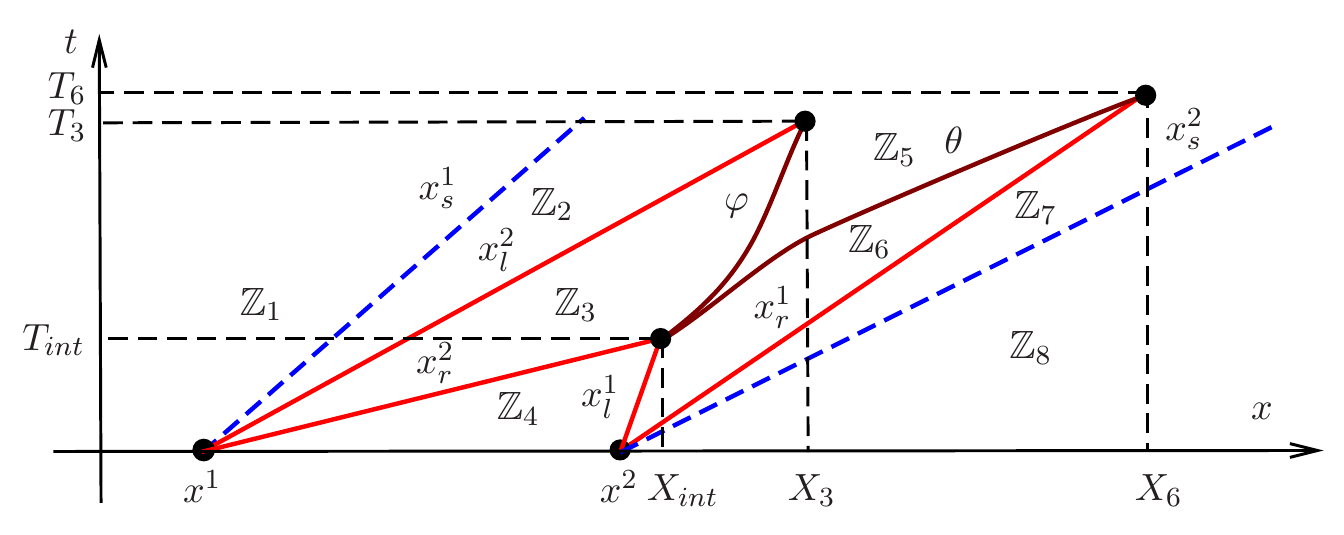}
  \caption{The zones in the $(x,t)$~--~plane. Weak discontinuities are solid lines, strong discontinuities~--- dashed lines.}
  \label{zhshel:fig1}
\end{figure}

For initial data  (\ref{zhshel:eq:3.02}) the evolution of zones on $(x,t)$-plane is shown in figure~\ref{zhshel:fig1}.

The fronts of the rarefaction waves and the shock waves are shown by solid and dashed lines, respectively. Figure~\ref{zhshel:fig1} corresponds to some finite time interval $(0, T_3)$. It is obvious that the  intersection of the boundaries corresponds to the interaction of discontinuities. There are only four types of interactions:
($\textbf{i}$) the interaction of two strong discontinuities,
($\textbf{ii}$, $\textbf{iii}$) the interaction of strong and weak discontinuities,
($\textbf{\textbf{iv}}$) the interaction of two weak discontinuities.
For instance, the point $(X_{int},T_{int})$ is interaction of two weak discontinuities $x_{r}^{2}$ and $x_{l}^{1}$. \label{interactions:1}

We recall that the waves arising at the breakup of initial discontinuity can be classified using the eigenvalues $\lambda^{k}$. The wave corresponding to the eigenvalue $\lambda^{k}$ is called k-wave. Thus,  the line $x_s^1$ corresponds to the 1-shock wave,
the lines $x_l^2$ and $x_r^2$  correspond to the left and right fronts of 2-rarefaction wave, etc. (see figure~\ref{zhshel:fig1}).

The solution of the problem (\ref{zhshel:eq:2.03}), (\ref{zhshel:eq:2.04}), (\ref{zhshel:eq:3.01})--(\ref{zhshel:eq:3.03}) in time interval $(0,T_{int})$ is described with the help of zones (see figure~\ref{zhshel:fig1})
\begin{eqnarray}\label{zhshel:eq:3.06}
\mathbb{Z}_1=\{\mu^1, \mu^2; -\infty, x_s^1(t)\},\quad
\mathbb{Z}_2=\{q^1, \mu^2; x_s^1(t), x_l^2(t)\},
\end{eqnarray}
\begin{eqnarray*}
\mathbb{Z}_3=\{q^1, R^2_a(z^2); x_l^2(t), x_r^2(t)\},\quad
\mathbb{Z}_4=\{q^1, q^2; x_r^2(t), x_l^1(t)\},
\end{eqnarray*}
\begin{eqnarray*}
\mathbb{Z}_6=\{R^1_a(z^1), q^2; x_l^1(t), x_r^1(t)\},\quad
\mathbb{Z}_7=\{\mu^1, q^2; x_r^1(t), x_s^2(t)\},
\end{eqnarray*}
\begin{eqnarray*}
\mathbb{Z}_8=\{\mu^1, \mu^2; x_s^2(t), +\infty\},
\end{eqnarray*}
where
\begin{eqnarray}\label{zhshel:eq:3.07}
R^1_a(z^1)=\sqrt{\frac{z^1}{q^2}}, \quad z^1=\frac{x-x^2}{t}, \quad R^2_a(z^2)=\sqrt{\frac{z^2}{q^1}}, \quad z^2=\frac{x-x^1}{t}.
\end{eqnarray}

The motion laws of the left and right rarefaction waves fronts (weak discontinuities) have the following form
\begin{eqnarray}\label{zhshel:eq:3.08}
x_l^1(t)=x^2 + q^1 q^1 q^2 t, \quad x_r^1(t)=x^2 + \mu^1 \mu^1 q^2 t,
\end{eqnarray}
\begin{eqnarray*}
x_l^2(t)=x^1 + q^1 \mu^2 \mu^2 t, \quad x_r^2(t)=x^1 + q^1 q^2 q^2 t.
\end{eqnarray*}

Strong discontinuities $x=x_s^1(t)$, $x=x_s^2(t)$ (shock waves) are moved with the velocities
$D^1=q^1 \mu^1 \mu^2$, $D^2=\mu^1 \mu^2 q^2$ as
\begin{eqnarray}\label{zhshel:eq:3.09}
x_s^1(t)=x^1 + q^1 \mu^1 \mu^2 t, \quad x_s^2(t)=x^2 + \mu^1 \mu^2 q^2 t.
\end{eqnarray}

As already mentioned, the solution (\ref{zhshel:eq:3.06})--(\ref{zhshel:eq:3.09}) of problem (\ref{zhshel:eq:2.03}), (\ref{zhshel:eq:2.04}),  (\ref{zhshel:eq:3.01})--(\ref{zhshel:eq:3.03}) is valid for the time interval $(0,T_{int})$. In the time moment $t=T_{int}$
the right front of 2-rarefaction wave $x=x_r^2(t)$ reaches the left front of 1-rarefaction wave $x=x_l^1(t)$. This moment corresponds to interaction of waves. The coordinates of the point $(X_{int},T_{int})$ on $(x,t)$-plane (see figure~\ref{zhshel:fig1}) can be calculated as the solution of the equation $x^2_r(t)=x^1_l(t)$
\begin{eqnarray*}
T_{int}=\frac{x^2-x^1}{q^1 q^2(q^2-q^1)}, \quad
X_{int}=\frac{x^1 q^1- x^2 q^2}{q^1 - q^2}.
\end{eqnarray*}

\begin{figure}[H]
  \centering
\includegraphics[scale=0.54]{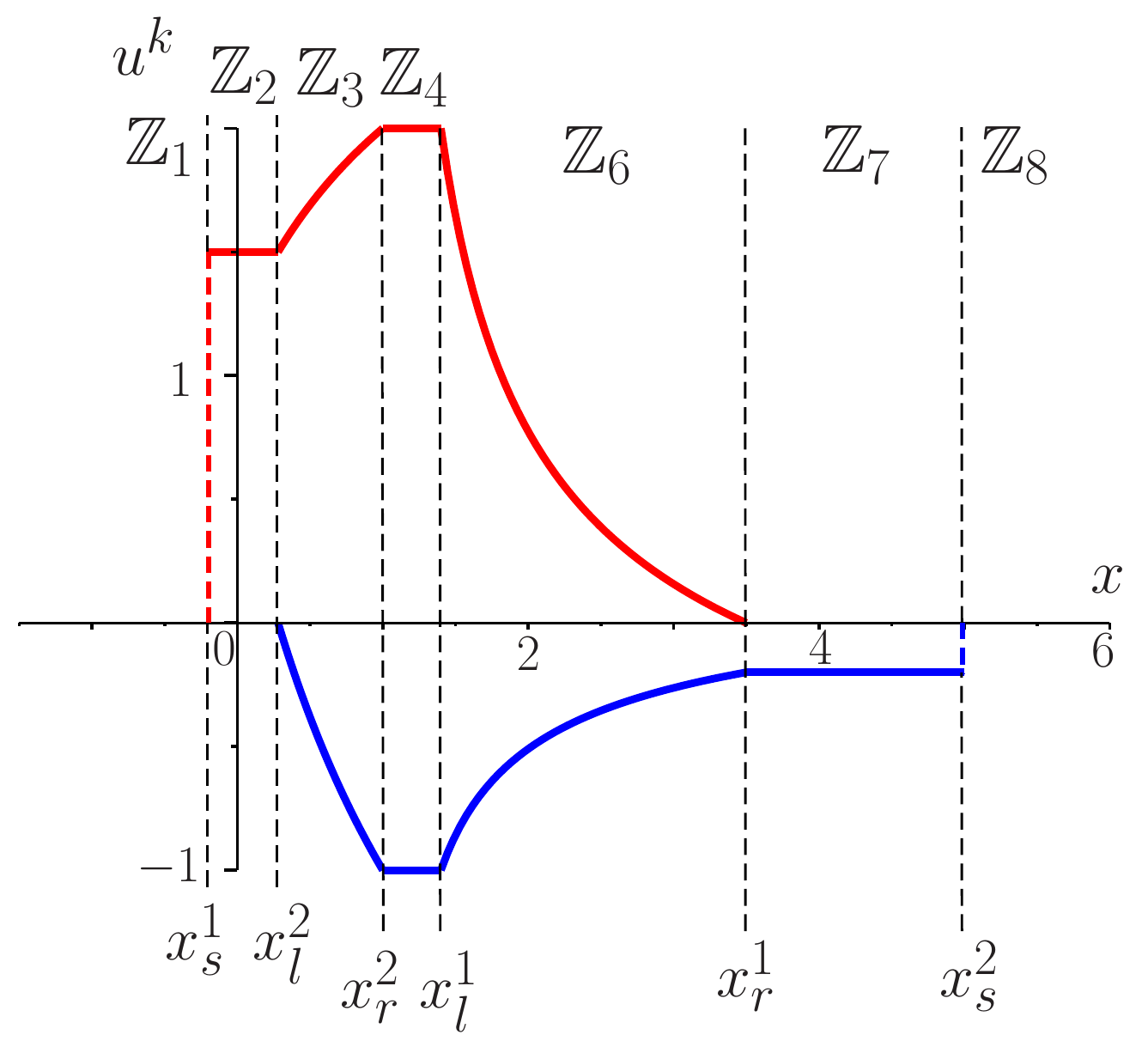}\qquad
\includegraphics[scale=0.54]{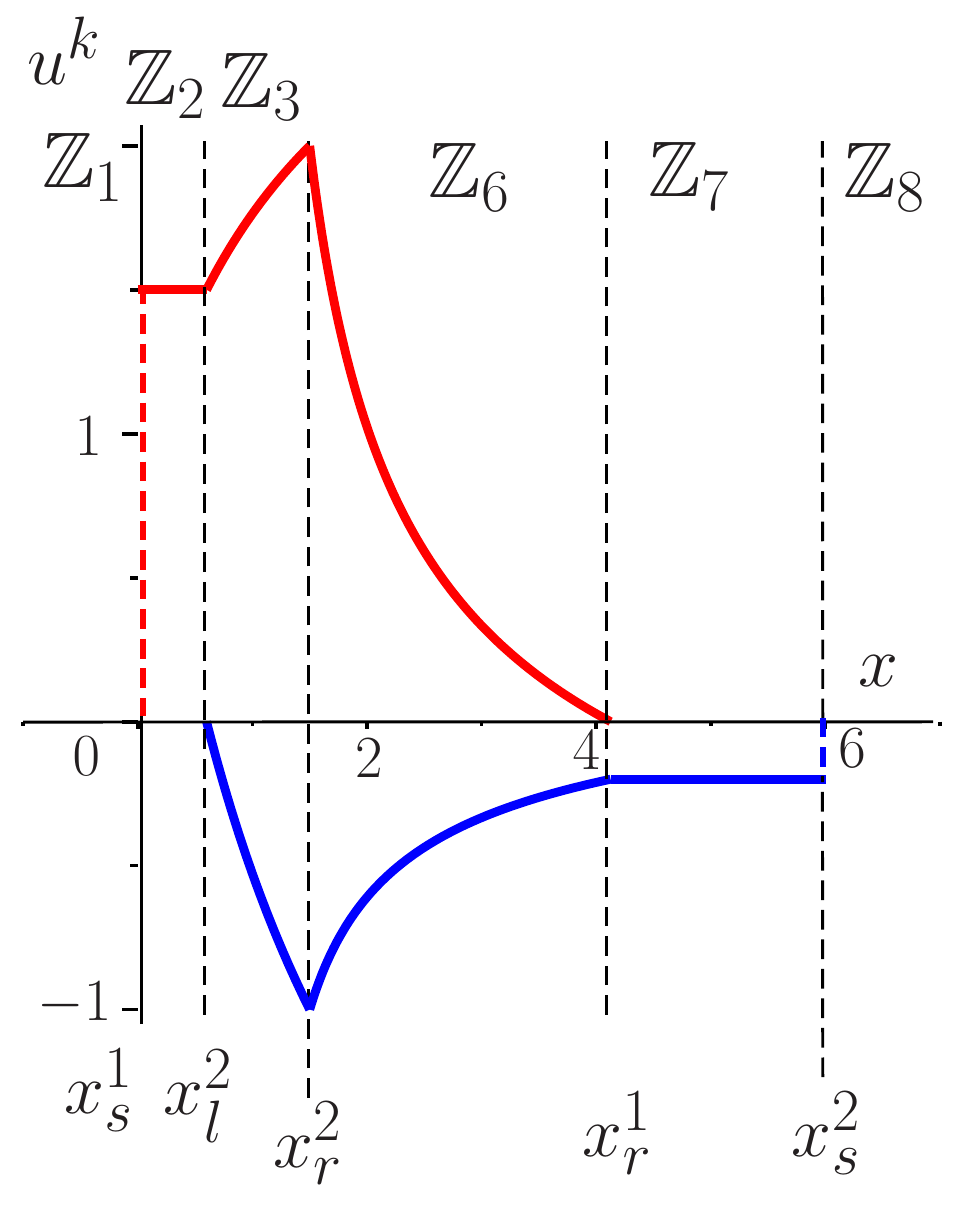}
  \caption{Concentrations $u^k$ at the moments $t_*=0.01$, $t_*=T_{int}=0.0125$}
  \label{zhshel:fig2}
\end{figure}

For completeness, the solutions of the original problem (\ref{zhshel:eq:2.01}) for concentrations $u^k$ at  moments  $t=T_{int}$ and $t=t_* < T_{int}$  are presented in figure~\ref{zhshel:fig2} in the case of the parameters
\begin{eqnarray}\label{zhshel:eq:3.10}
q^1=2, \quad \mu^1=5, \quad \mu^2=8, \quad q^2=10,\quad x^1=-1,\quad x^2=1.
\end{eqnarray}

As a result of interaction between two rarefaction waves the zone $\mathbb{Z}_4$ is disappeared. New zone $\mathbb{Z}_5$ is appeared
and the boundaries of the zones $\mathbb{Z}_3$, $\mathbb{Z}_6$ are modified (but the Riemann invariants in these zones are saved the same)
\begin{eqnarray*}
\mathbb{Z}_5=\{R^1(x,t), R^2(x,t); \varphi(t), \theta(t)\},
\end{eqnarray*}
\begin{eqnarray*}
\mathbb{Z}_3=\{q^1, R^2_a(z^2)); x_l^2(t), \varphi(t)\},\quad
\mathbb{Z}_6=\{R^1_a(z^2), q^2, \theta(t), x_r^1(t)\}.
\end{eqnarray*}
Here $\varphi(t)$, $\theta(t)$ are functions that determine new boundaries of the   zone.

Obviously, the lines $x=\varphi(t)$, $x=\theta(t)$ are weak discontinuities that are identified by corresponding characteristic directions $\lambda^1$ and $\lambda^2$. On these lines
the Riemann invariants $R^k$ have the weak discontinuity.

To determine, for example  $\varphi(t)$,  we have the Cauchy problem
\begin{eqnarray}\label{zhshel:eq:3.14}
\frac{d\varphi}{dt}=\lambda^1(q^1,R^2_a), \quad \varphi(T_{int})=X_{int}.
\end{eqnarray}
Integrating (\ref{zhshel:eq:3.14}) we get
\begin{eqnarray}\label{zhshel:eq:3.15}
\left(\varphi-x^1\right)^{1/2}=
\left(q^1\right)^{3/2}
\left(t^{1/2}-T_{int}^{1/2}\right)+
\left(X_{int}-x^1\right)^{1/2}.
\end{eqnarray}
Similarly, the line of weak discontinuity $x=\theta(t)$ is determined by formula
\begin{eqnarray}\label{zhshel:eq:3.16}
\left(\theta-x^2\right)^{1/2}=
\left(q^2\right)^{3/2}
\left(t^{1/2}-T_{int}^{1/2}\right)+
\left(X_{int}-x^2\right)^{1/2}.
\end{eqnarray}

Easily to describe the behavior of the lines $x=\varphi(t)$, $x=\theta(t)$ on the $(x,t)$-plane.
The zone $\mathbb{Z}_3$ (see figure~\ref{zhshel:fig1}) disappears at the moment $t=T_3$ in the point $x=X_3$ when characteristic
$x=\varphi(t)$ intersects characteristic $x=x_l^2(t)$.
The point $(X_3, T_3)$ is determined by  the relation $X_3=\varphi(T_3)=x_l^2(T_3)$
\begin{eqnarray*}
T_3=\frac{(q^2-q^1)^2}{(q^1-\mu^2)^2}\,T_{int}, \quad X_3=x^1+q^1 \mu^2 \mu^2 T_3.
\end{eqnarray*}

The zone $\mathbb{Z}_6$ (see figure~\ref{zhshel:fig1}) disappears at the moment $t=T_6$ in the point $x=X_6$, when characteristic
$x=\theta(t)$  intersects characteristic $x=x_r^1(t)$. Point $(X_6, T_6)$ is determined by the relation $X_6=\theta(T_6)=x_r^1(T_6)$
\begin{eqnarray*}
T_6=\frac{(q^2-q^1)^2}{(q^2-\mu^1)^2}\,T_{int}, \quad X_6=x^2+ \mu^1 \mu^1 q^2  T_6.
\end{eqnarray*}

The results listed in this section are given in one form or another in the papers \cite{ElaevaIzvestiya,ElaevaMM,Elaeva_Diss,Elaeva_ZhVM}.
In particular, these papers indicate that further study of interactions between waves requires the solution of the Goursat problem with discontinuous initial data on characteristics. Primarily, it concerns to the interaction between weak discontinuities.

An effort to solve a similar problem with the help of the numerical methods for the modified diffusive approximation has shown that the numerical results do not allow us to describe adequately the interaction of weak discontinuities (see \cite{Elaeva_Diss}).

In the following sections we show that the Goursat problem can be solved using a variant of the hodograph method based on the conservation laws (see \cite{SenashovYakhno} and also \cite{Zhuk_Shir_ArXiv_2014_Part2}). Moreover, this method allows us to study any interaction of the discontinuities  for the original problem (\ref{zhshel:eq:2.03}), (\ref{zhshel:eq:2.04}), (\ref{zhshel:eq:3.01}), (\ref{zhshel:eq:3.02}) in the case when inequalities (\ref{zhshel:eq:3.03}) are valid.

\subsection{The Goursat problem}\label{zhshel:sec:3.1}

To determine the Riemann invariant $R^k(x,t)$ in zone $\mathbb{Z}_5$ we solve the Goursat problem for the equations
(\ref{zhshel:eq:2.03}), (\ref{zhshel:eq:2.04}) with data that have weak discontinuities on the characteristics $x=\varphi(t)$, $x=\theta(t)$
\begin{eqnarray}\label{zhshel:eq:3.19}
R^1\bigr|_{x=\varphi(t)}=q^1, \quad R^2\bigr|_{x=\varphi(t)}=R^2_a(z^2(\varphi(t),t)),
\end{eqnarray}
\begin{eqnarray}\label{zhshel:eq:3.20}
R^1\bigr|_{x=\theta(t)}=R^1_a(z^1(\theta(t),t)), \quad R^2\bigr|_{x=\theta(t)}=q^2.
\end{eqnarray}
Note that  data (\ref{zhshel:eq:3.19}), (\ref{zhshel:eq:3.20}) are matched at the point
$(X_{int},T_{int})$
\begin{eqnarray}\label{zhshel:eq:3.21}
R^1\bigr|_{t=T_{int}}=q^1, \quad R^2\bigr|_{t=T_{int}}=q^2.
\end{eqnarray}

To avoid confusions we note that the problem (\ref{zhshel:eq:2.03}), (\ref{zhshel:eq:2.04}), (\ref{zhshel:eq:3.19})--(\ref{zhshel:eq:3.21}) is
the initial-boundary value problem for hyperbolic equations. Analysis of the solution methods of such problems one can find, for example,
in~\cite{RozhdestvenskiiYanenko}.
This problem has some originality in view of the fact that the boundary conditions (\ref{zhshel:eq:3.19}), (\ref{zhshel:eq:3.20}) are specified on characteristics,
but the initial condition (\ref{zhshel:eq:3.21}) is specified only in one point $(X_{int},T_{int})$. 



\section{Hodograph method}\label{zhshel:sec:4}

We recall that the hodograph method allows us to transform a quasi-linear system of equation in a linear system, using the replacement  of dependent and independent variables:  $(R^1,R^2) \rightleftharpoons (x,t)$. To determine $t(R^1,R^2)$, $x(R^1,R^2)$ we have the linear systems of equations
\begin{eqnarray}\label{zhshel:eq:4.01}
x_{R^2}- \lambda^1 t_{R^2}=0, \quad  x_{R^1}- \lambda^2 t_{R^1}=0.
\end{eqnarray}
The compatibility condition of these equations leads to a linear hyperbolic system with variable coefficients.
Taking into account (\ref{zhshel:eq:2.04}), for the function $t(R^1,R^2)$ we have
\begin{eqnarray}\label{zhshel:eq:4.02}
t_{R^1 R^2} + \frac{2}{R^2- R^1} (t_{R^1}-t_{R^2})=0.
\end{eqnarray}
The equation (\ref{zhshel:eq:4.02}) has the Riemann-Green function (see, e.g. \cite{Copson})
\begin{eqnarray*}
V(r^1,r^2|R^1,R^2)=
\frac{((R^1 + R^2)(r^1 + r^2)-2(R^1R^2 + r^1r^2))(r^1 - r^2)}{(R^1 - R^2)^3},
\end{eqnarray*}
which satisfies to equation (\ref{zhshel:eq:4.02}) for the variables $R^1$, $R^2$ and
satisfies to conjugate equation for the variables $r^1$, $r^2$.
In addition the condition $V(r^1,r^2|r^1,r^2)=1$ is valid.

The explicit form of the Riemann-Green function allows us to obtain the solution of the Goursat
problem for equation (\ref{zhshel:eq:4.02}) as (see, e.g.
\cite[formula (2.20)]{Bizadze})
\begin{eqnarray}\label{zhshel:eq:4.04}
t(R^1,R^2)=
-V(R^1_0,R^2_0,R^1,R^2)t(R^1_0, R^2_0)+{}
\end{eqnarray}
\begin{eqnarray*}
{}+\frac{2(R^1-R^1_0)(R^2-R^1_0)}{(R^1-R^2)^3}H^2(R^2)
-\frac{2(R^1-R^2_0)(R^2-R^2_0)}{(R^1-R^2)^3}H^1(R^1)
\end{eqnarray*}
\begin{eqnarray*}
{}+\frac{(R^1-R^2_0)^2}{(R^1-R^2)^2}
t(R^1,R^2_0)
+\frac{(R^2-R^1_0)^2}{(R^1-R^2)^2}
t(R^1_0,R^2),
\end{eqnarray*}
where
\begin{eqnarray}\label{zhshel:eq:4.05}
H^1(R^1)=\int\limits_{R^1_0}^{R^1}t(\tau,R^2_0)\,d\tau,\quad
H^2(R^2)=\int\limits_{R^2_0}^{R^2}t(R^1_0,\tau)\,d\tau.
\end{eqnarray}
Here $t(R^1,R^2_0)$, $t(R^1_0,R^2)$ are the function $t(R^1,R^2)$ on characteristics $R^2=R^2_0$ and $R^1=R^1_0$, that one can get using  conditions (\ref{zhshel:eq:3.19}), (\ref{zhshel:eq:3.20}).

The condition (\ref{zhshel:eq:3.19}) means that for characteristic $x=\varphi(t)$ the Riemann invariant is $R^1=q^1$.
Using the relations (\ref{zhshel:eq:3.19}), (\ref{zhshel:eq:3.07}), (\ref{zhshel:eq:3.15})
and excluding function $\varphi(t)$  we get
\begin{eqnarray}\label{zhshel:eq:4.06}
t=\frac{(q^2 - q^1)^2}{(R^2 - q^1)^2}T_{int}=t(R^1_0,R^2), \quad R^1_0=q^1.
\end{eqnarray}
Similarly, the condition (\ref{zhshel:eq:3.20}) leads to (it is enough to replace indices ${1 \rightleftharpoons 2}$)
\begin{eqnarray}\label{zhshel:eq:4.07}
t=\frac{(q^1 - q^2)^2}{(R^1 - q^2)^2}T_{int}=t(R^1,R^2_0), \quad R^2_0=q^2.
\end{eqnarray}
Finaly, substituting (\ref{zhshel:eq:4.06}), (\ref{zhshel:eq:4.07}) in~(\ref{zhshel:eq:4.04}), (\ref{zhshel:eq:4.05}) we obtain
\begin{eqnarray}\label{zhshel:eq:4.08}
t(R^1,R^2)=T_{int}V(q^1,q^2|R^1,R^2)=
\end{eqnarray}
\begin{eqnarray*}
=\frac{(x^2-x^1)(2 R^1 R^2 + 2 q^1 q^2 - (q^1 +  q^2)(R^1 + R^2))}{q^1 q^2 (R^1-R^2)^3}.
\end{eqnarray*}
Note that the initial condition (\ref{zhshel:eq:3.21}) is automatically satisfied because the Riemann-Green function has properties:
$t(q^1,q^2)=T_{int}V(q^1,q^2|q^1,q^2)=T_{int}$.

To obtain the function $x(R^1,R^2)$ it is enough to pay attention that after the substitutions $R^i=1/K^i$ (see \cite{Pavlov_Preprint})
the equation for $x(R^1,R^2)$
\begin{eqnarray*}
x_{K^1 K^2} + \frac{2}{K^2- K^1} (x_{K^1}-x_{K^2})=0
\end{eqnarray*}
formally coincides with the equation (\ref{zhshel:eq:4.02}).

Performing the replacements $t \rightleftharpoons x$, $R^i \rightleftharpoons 1/K^i$  in the formula (\ref{zhshel:eq:4.08}) leads to
\begin{eqnarray}\label{zhshel:eq:4.11}
x(R^1,R^2)=
(x^2-x^1)(R^1 R^2)^2 \frac{R^2  +  R^1 - 2  ( q^1+q^2 ) }{q^1 q^2 (R^1-R^2)^3}+{}
\end{eqnarray}
\begin{eqnarray*}
{}+\frac{(x^1 (R^1)^3 - x^2 (R^2)^3 )  + 3  R^1 R^2 ( R^2 x^2 - R^1 x^1)  }{(R^1-R^2)^3}.
\end{eqnarray*}

Thus, the relations (\ref{zhshel:eq:4.08}), (\ref{zhshel:eq:4.11}) give us the implicit  solution of the Goursat problem
(\ref{zhshel:eq:2.03}), (\ref{zhshel:eq:2.04}), (\ref{zhshel:eq:3.19})--(\ref{zhshel:eq:3.21}).


\section{Construction of the solution on isochrone}\label{zhshel:sec:5}

The functions $R^1(x,t)$, $R^2(x,t)$ that are inverse to the functions $t(R^1,R^2)$, $x(R^1,R^2)$ allow us to obtain the Riemann invariants for the zone $\mathbb{Z}_5$.
Relations  describe the solution to $R^1(x,t)$, $R^2(x,t)$ on the interval $T_{int} \leqslant t \leqslant \min(T_3,T_6)$ only. This should be kept in mind because, from the moment $T_{int}$, the zone $\mathbb{Z}_5$ disappears as a result of weak discontinuities interaction. We also note that the direct solution of algebraic equations $t=t(R^1,R^2)$, $x=x(R^1,R^2)$, in practice, is much more difficult than the solution of differential equations which are presented below.

To obtain the differential equations for functions $R^1(x,t)$, $R^2(x,t)$ we differentiate functions $t=t(R^1,R^2)$, $x=x(R^1,R^2)$ with respect to $x$
\begin{eqnarray*}
0=  t_{R^1} R^1_x +  t_{R^2} R^2_x, \quad 1=  x_{R^1}R^1_x +  x_{R^2}R^2_x.
\end{eqnarray*}
Obviously, these equations are valid on any level line $t_*=t(R^1,R^2)$ (isochrone), where $t_*$ is fixed. For isochrone this system is the ODEs which can be  written in the following form
\begin{eqnarray}\label{zhshel:eq:5.01}
\frac{d R^1(x,t_*)}{dx}=-\frac{t_{R^2}(R^1,R^2)}{\Delta(R^1,R^2)},\quad
\frac{d R^2(x,t_*)}{dx}=\frac{t_{R^1}(R^1,R^2)}{\Delta(R^1,R^2)},
\end{eqnarray}
where
\begin{eqnarray}\label{zhshel:eq:5.02}
\Delta=t_{R^1}x_{R^2}-t_{R^2}x_{R^1}=(\lambda^1-\lambda^2)t_{R^1}t_{R^2}.
\end{eqnarray}
Here, to exclude $x_{R^1}$, $x_{R^1}$, we use (\ref{zhshel:eq:4.01}).

For the equations (\ref{zhshel:eq:5.01}), (\ref{zhshel:eq:5.02}) the initial conditions result from the relations (\ref{zhshel:eq:3.19})
\begin{eqnarray}\label{zhshel:eq:5.03}
R^1(x_*,t_*)=q^1, \quad R^2(x_*,t_*)=R^2_a(x_*,t_*),\quad
x_*=\varphi(t_*).
\end{eqnarray}
Integrating the Cauchy problem (\ref{zhshel:eq:5.01})--(\ref{zhshel:eq:5.03}) from $\varphi(t_*)$ to $\theta(t_*)$ we get the solution
$R^1(x,t_*)$, $R^2(x,t_*)$, $\varphi(t_*)  \le x \le  \theta(t_*)$ for any $t_*$.

We note that the derivatives $t_{R^1}$, $t_{R^2}$ in the right part of equations (\ref{zhshel:eq:5.01}), (\ref{zhshel:eq:5.02}) can be calculated in the explicit form with the help (\ref{zhshel:eq:4.08}) only. The formula (\ref{zhshel:eq:4.11}) for function $x(R^1,R^2)$  is not \textbf{required}!

\begin{figure}[H]
\centering
\includegraphics[scale=0.48]{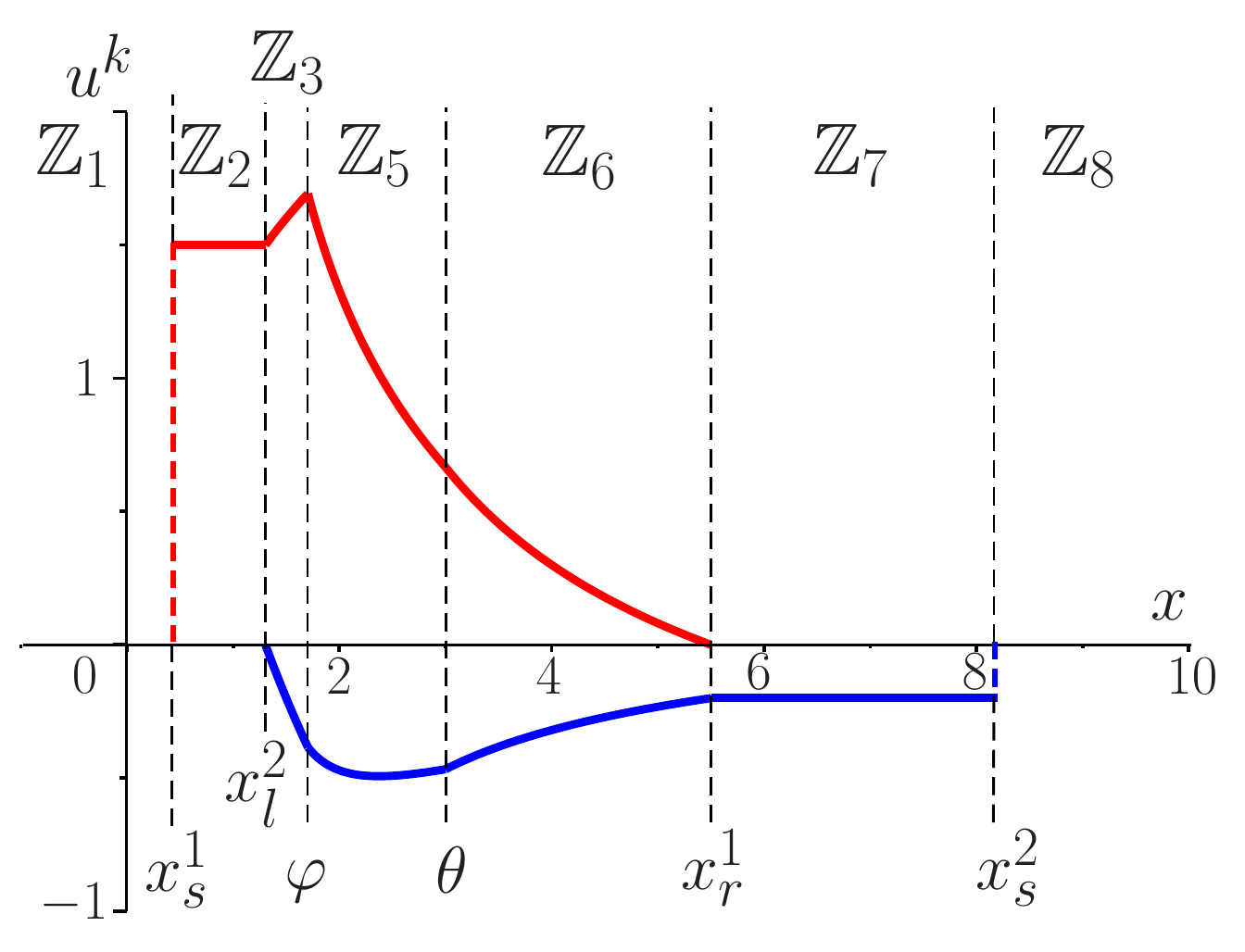}
\includegraphics[scale=0.48]{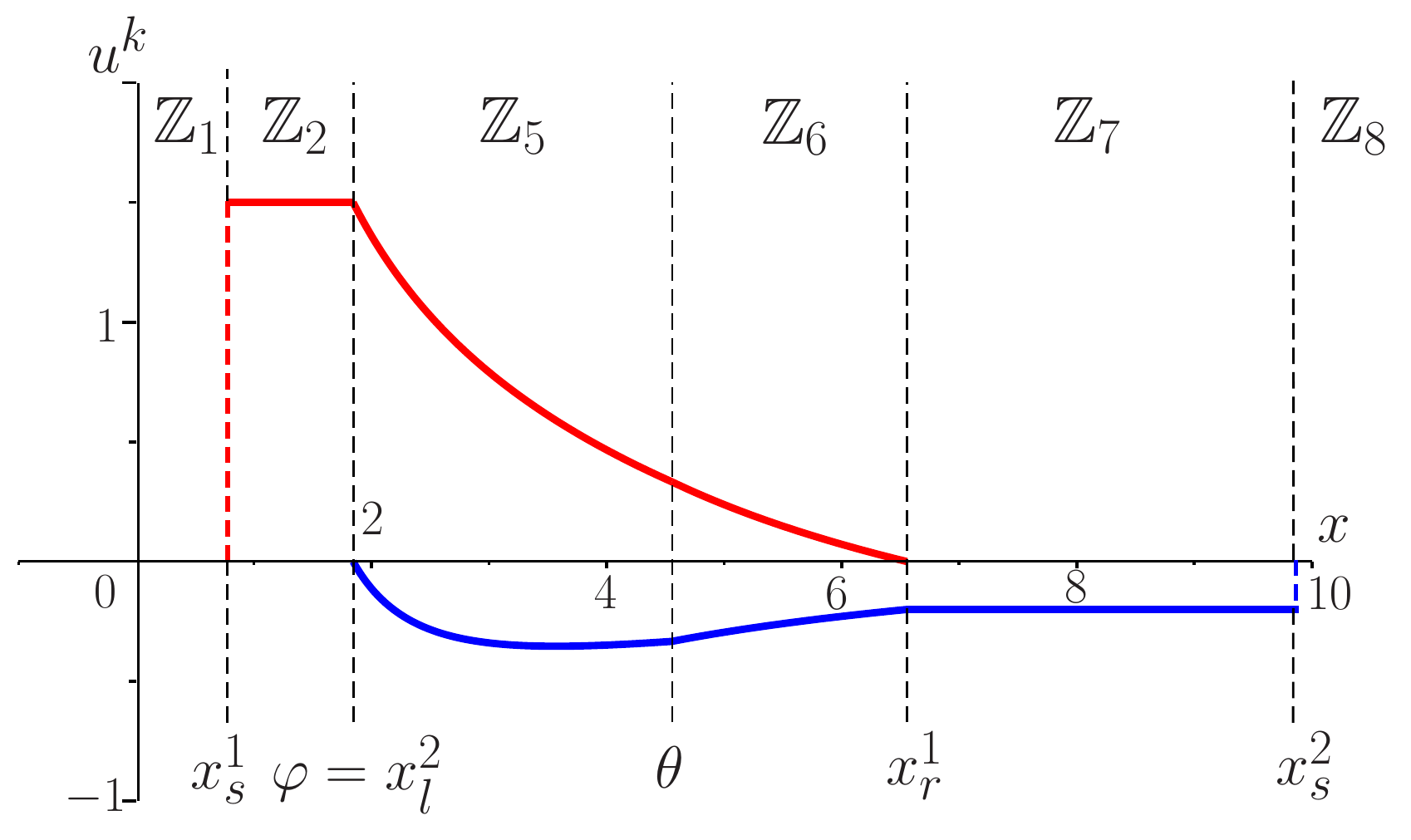}
\caption{Concentrations $u^k$ at the moments $t_*=0.018$, $t_*=T_3=0.022\dots=1/45$}.
\label{zhshel:fig3}
\end{figure}

For parameters (\ref{zhshel:eq:3.10}) the numerical results for the Cauchy problem (\ref{zhshel:eq:5.01})--(\ref{zhshel:eq:5.03}) are presented
in figure~\ref{zhshel:fig3}


\section{Interaction of discontinuities at $t \geqslant T_3$ (and at $t \geqslant T_6$ )}\label{zhshel:sec:6}

\begin{figure}[H]
  \centering
\includegraphics[scale=0.90]{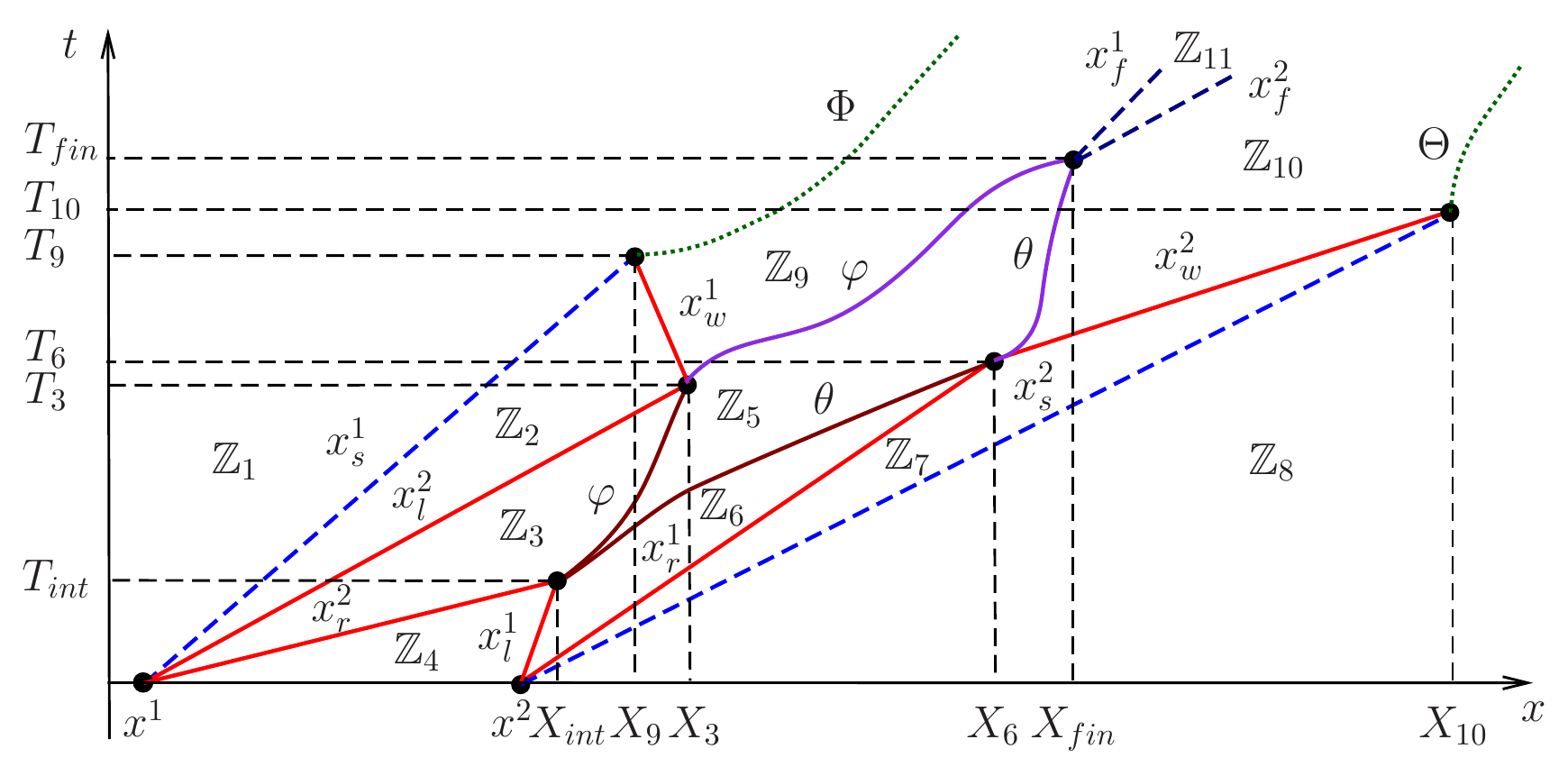}
  \caption{Evolution of initial discontinuities}
  \label{zhshel:fig4}
\end{figure}

Complete scenario  of the solution behavior is shown on the figure~\ref{zhshel:fig4}.
To investigate the behavior of the zone boundaries after the interaction it should be remembered that the weak discontinuities of the solution do not disappear over time and move along characteristics.
It means that the boundaries $x=x_w^1(t)$  and  $x=\varphi(t)$ of the new zone $\mathbb{Z}_9$ are the characteristics that have the weak discontinuity of the Riemann invariants $R^1$ and $R^2$ correspondingly. Similarly, the boundaries $x=\theta(t)$ and $x=x_w^2(t)$ of the new zone $\mathbb{Z}_{10}$ are the characteristics that have the weak discontinuity of the Riemann invariants $R^1$ and $R^2$. Note that the line of weak discontinuity for the  Riemann invariant $R^1$ is denoted as $x=x^1_w(t)$, and the line of weak discontinuity for the  Riemann invariant $R^2$ is denoted as $x=x^2_w(t)$. The lines of discontinuities for the Riemann invariants $R^1$ and $R^2$ are still denoted as $x=\theta(t)$ and $x=\varphi(t)$. Obviously, for zones $\mathbb{Z}_{9}$, $\mathbb{Z}_{10}$  the functions $\theta(t)$,  $\varphi(t)$ are the new functions that do not coincide with previous.

We focus to the solving of the problem only in zone $\mathbb{Z}_{9}$. For zone $\mathbb{Z}_{10}$ the solution is obtained in a similar way.
Zone $\mathbb{Z}_{9}$ arises at time $t=T_3$, when the line $x=x_l^2(t)$ ($R^2$ weak discontinuity) intersects with line $x=\varphi(t)$ (also $R^2$ weak discontinuity).
In the moment of interaction the weak discontinuity of the Reamann invariant $R^1$ is appeared on line $x=x_w^1(t)$ and the weak discontinuity of the Reamann invariant $R^2$ is saved on line $x=\varphi(t)$.

Note that the invariant $R^2$ in the zone $\mathbb{Z}_{9}$ has the same value as in zone $\mathbb{Z}_{2}$, because line $x=x_w^1(t)$ is discontinuity line only for  invariant $R^1$. Hence, for zone $\mathbb{Z}_{9}$ we have $R^2=\mu^2$ (similarly for zone $\mathbb{Z}_{10}$ we have $R^1=\mu^1$).
Thus, in the zone $\mathbb{Z}_{9}$ varies the invariant $R^1$ only and one can  consider the equation
\begin{eqnarray*}
R^1_t +\lambda^1(R^1,\mu_2)R^1_x=0
\end{eqnarray*}
only for  invariant $R^1$ instead of the system.

Easily to formulate the continuity conditions for invariant $R^1$ (and $R^2$) on the lines $x=x_w^1(t)$, $x=\varphi(t)$
\begin{eqnarray}\label{zhshel:eq:6.01}
R^1\bigr|_{x=x^1_w(t)}=q^1, \quad R^2\bigr|_{x=x^1_w(t)}=\mu^2,
\end{eqnarray}
\begin{eqnarray}\label{zhshel:eq:6.02}
R^1\bigr|_{x=\varphi(t)}=R^1_5\bigr|_{x=\varphi(t)}, \quad
\mu^2=R^2_5\bigr|_{x=\varphi(t)},
\end{eqnarray}
where $R^1_5(x,t)$,  $R^2_5(x,t)$ are the invariants for zone $\mathbb{Z}_{5}$ that are obtained in section~\ref{zhshel:sec:5}.

The conditions (\ref{zhshel:eq:6.01}) mean that on the line $x=x^1_w(t)$ we have $\lambda^1(R^1,R^2)=q^1 q^1 \mu^2$.  It allows us to get the function $x^1_w(t)$ (characteristic), i.e the motion law for left boundary of zone $\mathbb{Z}_{9}$,
\begin{eqnarray*}
x^1_w(t)=X_3+ q^1 q^1 \mu^2 (t-T_3).
\end{eqnarray*}

The right boundary of zone $\mathbb{Z}_{9}$, i.e. function $\varphi$, can be represented  in parametric form with the help of
the explicit relations (\ref{zhshel:eq:4.08}), (\ref{zhshel:eq:4.11}) for $t(R^1,R^2)$ and $x(R^1,R^2)$ and the continuity conditions (\ref{zhshel:eq:6.02}) for the Riemann invariant $R^2$
\begin{eqnarray}\label{zhshel:eq:6.04}
\varphi(t)=x(\rho^1,\mu^2), \quad t=t(\rho^1,\mu^2),  \quad q^1 \leqslant \rho^1 \leqslant \mu^1,
\end{eqnarray}
where $\rho^1$ is the parameter (function $R^1$ on the line $x=\varphi(t)$, i.e. $R^1_5(\varphi(t),t)$, plays the role of parameter).

Similarly, parametric representation of the function $\theta(t)$ is
\begin{eqnarray}\label{zhshel:eq:6.05}
\theta(t)=x(\mu^1,\rho^2), \quad t=t(\mu^1,\rho^2), \quad  \mu^2 \leqslant \rho^2 \leqslant q^2,
\end{eqnarray}
where $\rho^2$ is the parameter.

To obtain the Riemann invariant $R^1(x,t)$ for zone $\mathbb{Z}_9$ we have the Cauchy problem (classical method of characteristic)
\begin{eqnarray}\label{zhshel:eq:6.06}
\frac{dR^1}{dt}=0, \quad \frac{dx}{dt}=\lambda^1(R^1,\mu^2), \quad \frac{d}{dt}=\frac{\partial}{\partial t}+\lambda^1(R^1,\mu^2)\frac{\partial}{\partial x},
\end{eqnarray}
\begin{eqnarray*}
R^1\bigr|_{t=\tau}=R^1_5(\varphi(\tau),\tau), \quad x\bigr|_{t=\tau}=\varphi(\tau).
\end{eqnarray*}
Obviously, that the solution of problem (\ref{zhshel:eq:6.06}) is
\begin{eqnarray}\label{zhshel:eq:6.07}
R^1(x,t)=R^1_5(\varphi(\tau),\tau), \quad \tau=\tau(x,t),
\end{eqnarray}
where $\tau(x,t)$ is implicitly determined  by the relation
\begin{eqnarray}\label{zhshel:eq:6.08}
x=\varphi(\tau) + R^1_5(\varphi(\tau),\tau) R^1_5(\varphi(\tau),\tau)\mu^2 (t-\tau).
\end{eqnarray}

Thus, the zone $\mathbb{Z}_9$ is completely determined by formula
\begin{eqnarray*}
\mathbb{Z}_9=\{(R^1_5(\varphi(\tau),\tau), \mu^2), (x^1_w(t),\varphi(t))\},
\end{eqnarray*}
where for any point $(x,t)$ parameter $\tau$ is given by  (\ref{zhshel:eq:6.08}).

\begin{figure}[H]
  \centering
\includegraphics[scale=0.65]{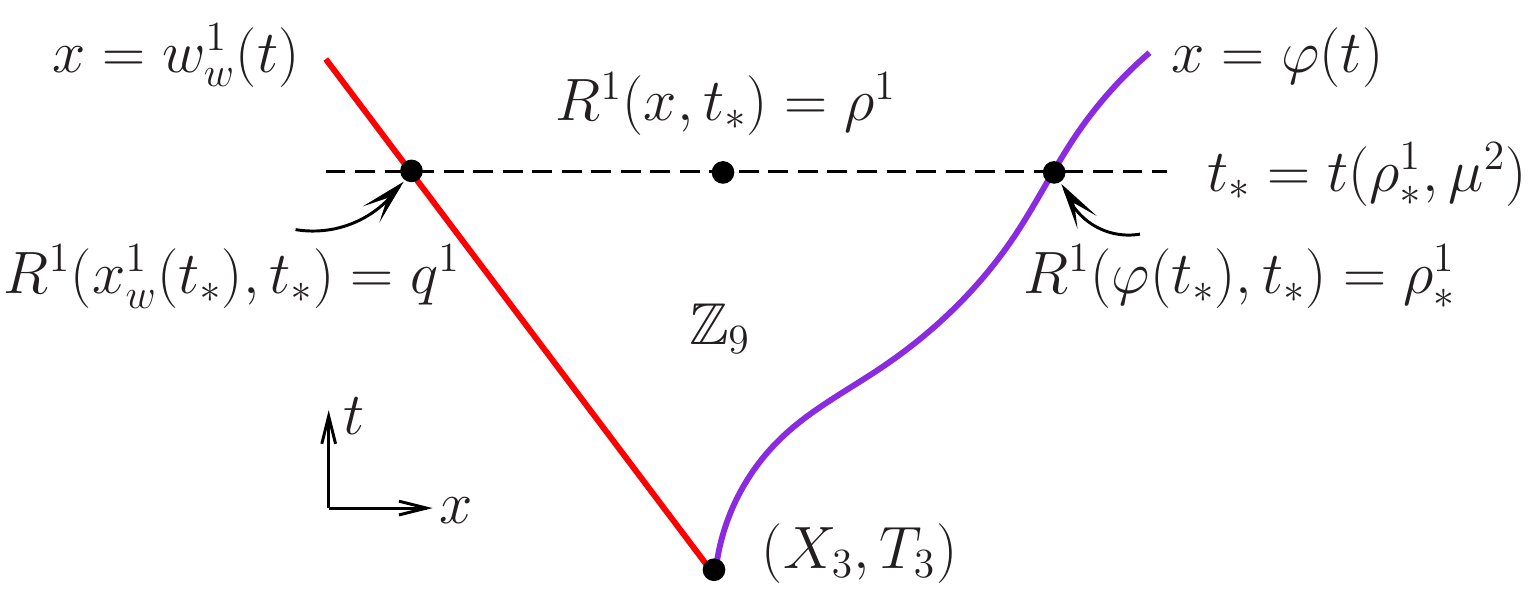}
  \caption{Fragment of the zone $\mathbb{Z}_9$}
  \label{zhshel:fig5}
\end{figure}

At first sight the relations obtained are complicated. In actuality, using again the idea of isochrone, it is easy to construct a solution which depends on only one parameter $\rho^1$ (see figure~\ref{zhshel:fig5}).

We fix a time $t_*$ and consider the isochrone $t_*=t(R^1,R^2)$. For zone $\mathbb{Z}_9$ we have $t_*=t(R^1,\mu^2)$. We introduce notation for the Riemann invariants on the isochrone: $R^1(x,t_*)=\rho^1$. Obviously (see figure~\ref{zhshel:fig5}), that the parameter $\rho^1$ satisfies to inequality $q^1 \leqslant \rho^1 \leqslant \rho^1_*$, where $\rho^1_*$ is root of equation
\begin{eqnarray}\label{zhshel:eq:6.09}
t_*=t(\rho^1_*,\mu^2), \quad q^1 \leqslant \rho^1_* \leqslant \mu^1.
\end{eqnarray}
This root $\rho^1_*$ is the value of Riemann invariant which corresponds to the intersection of isochrone $t_*=t(R^1,\mu^2)$ and  line $x=\varphi(t)$ at moment $t_*$.
The root $\rho^1_*$ always exists, because the equation (\ref{zhshel:eq:6.09}) coincides with one equation of system (\ref{zhshel:eq:6.04}), that determine the  function $\varphi(t)$. It is easy to check that the equation (\ref{zhshel:eq:6.09}) is a cubic equation. The solution of this equation should be chosen from the interval
$[q^1, \mu^1]$.

After replacement  $R^1_5(\varphi(\tau),\tau) \to \rho^1$ in the relations (\ref{zhshel:eq:6.07}), (\ref{zhshel:eq:6.08}) we get the one-parametric representation of the Riemann invariant  $R^1$ for any point $(\overline{x},t_*)$
\begin{eqnarray}\label{zhshel:eq:6.10}
R^1(\overline{x},t_*)=\rho^1, \quad q^1 \leqslant \rho^1 \leqslant \rho^1_* \leqslant \mu^1,
\end{eqnarray}
where
\begin{eqnarray}\label{zhshel:eq:6.11}
\overline{x}=\varphi(\tau) + \rho^1 \rho^1 \mu^2 (t_*-\tau), \quad \tau=t(\rho^1, \mu^2), \quad \varphi(\tau)=x(\rho^1, \mu^2).
\end{eqnarray}
It is easy to check that the first condition (\ref{zhshel:eq:6.01}) is satisfied at $\rho^1=q^1$.

Similarly, one can obtain solution for zone $\mathbb{Z}_{10}$. To write this solution it is enough  to replace
the indices  $1 \leftrightarrows 2$ and  the functions $\varphi \leftrightarrows \theta$.

\begin{figure}[H]
\centering
\includegraphics[scale=0.54]{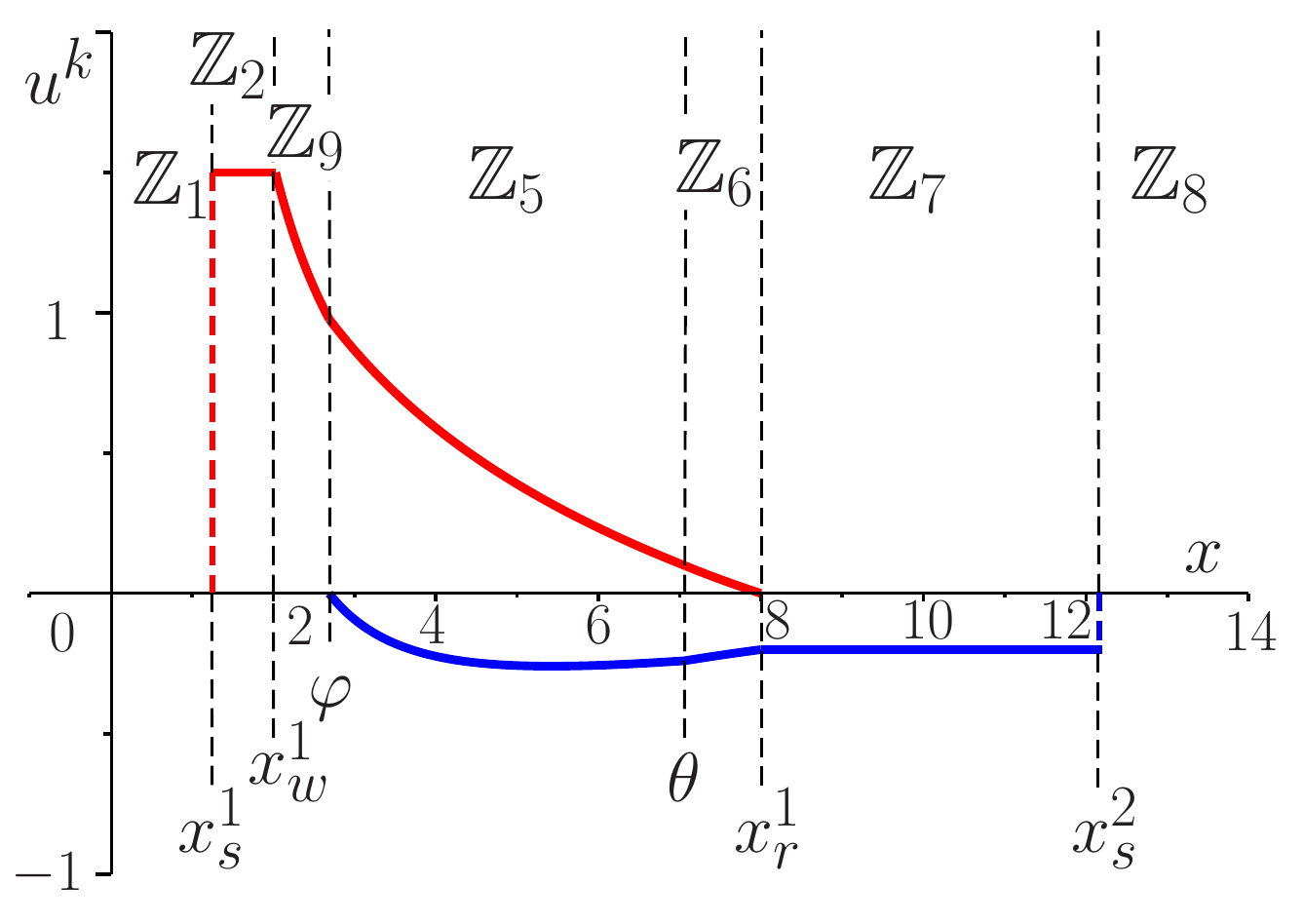}
\includegraphics[scale=0.54]{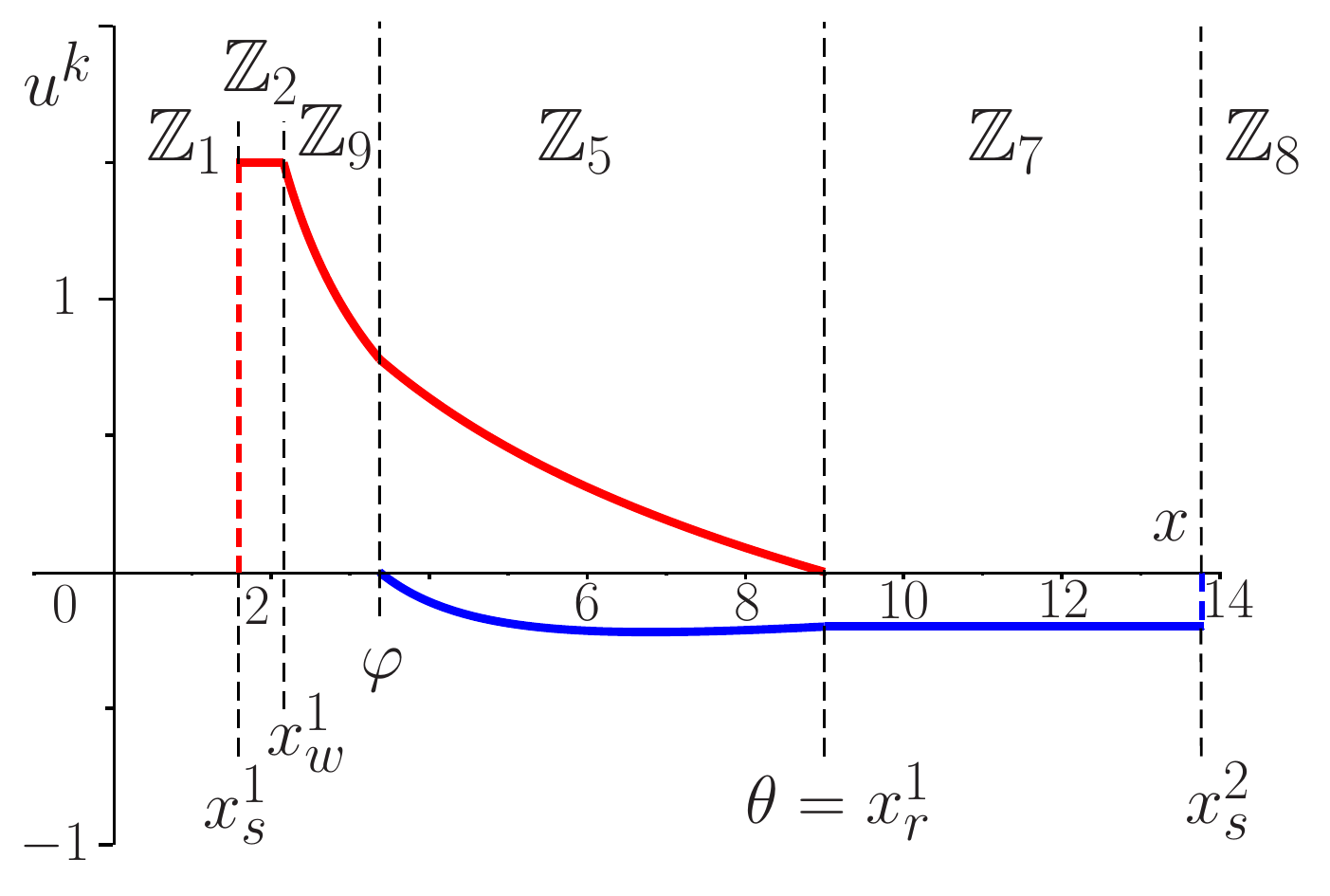}
\caption{Concentrations $u^k$ at the moments $t_*=0.028$, $t_*=T_6=0.032$}
\label{zhshel:fig6}
\end{figure}

\begin{figure}[H]
\centering
\includegraphics[scale=0.54]{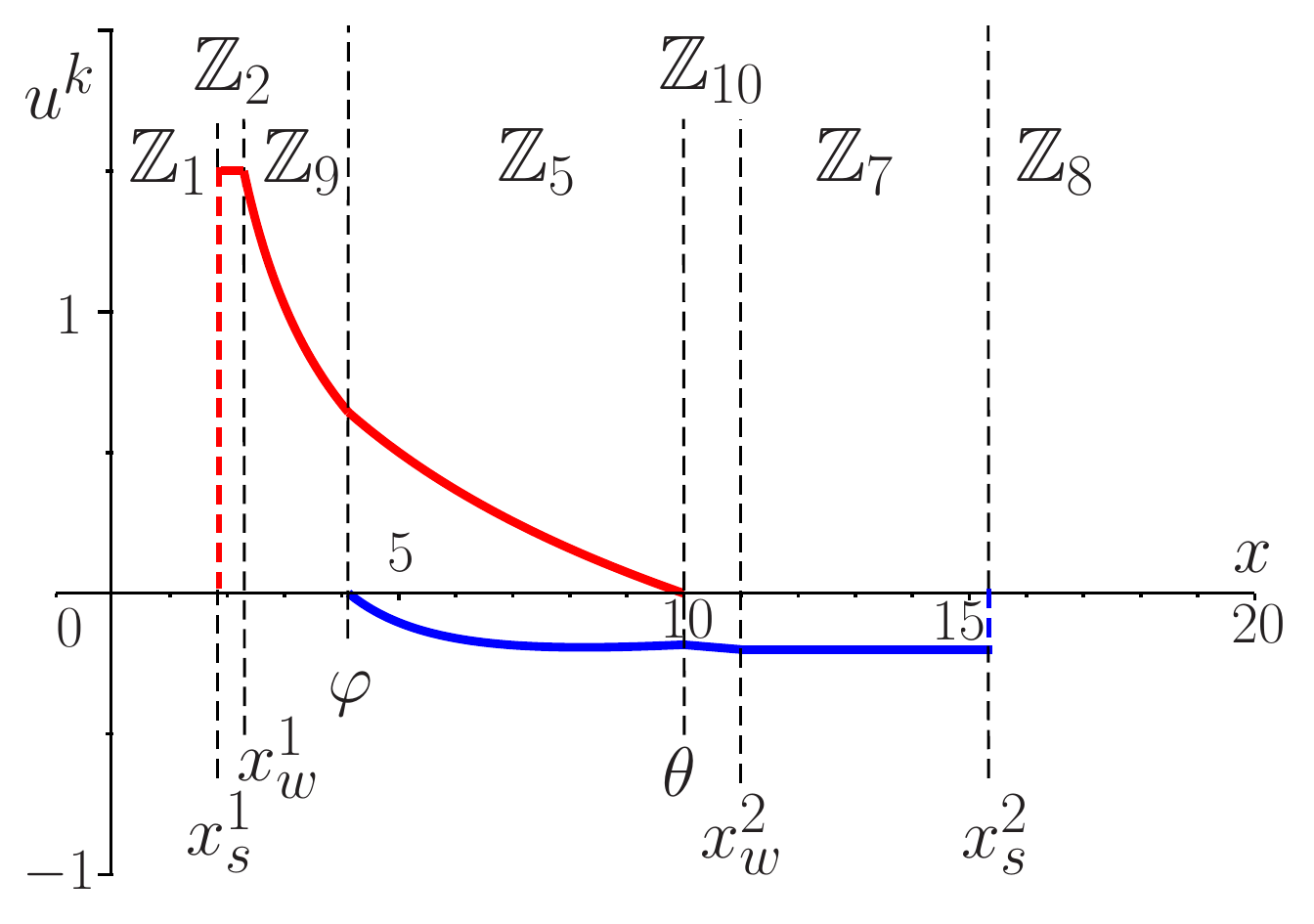}
\includegraphics[scale=0.54]{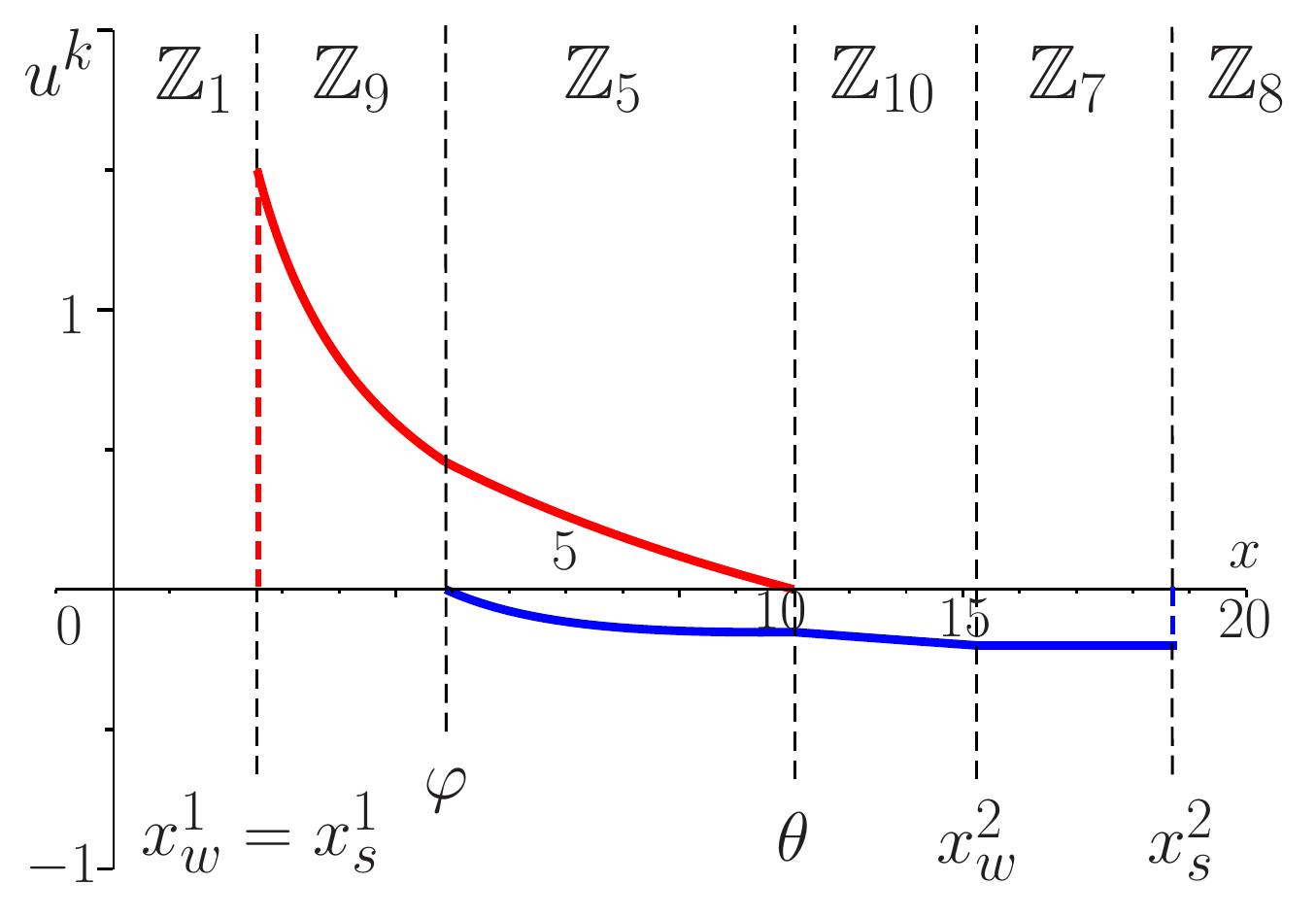}
\caption{Concentrations $u^k$ at the moments $t_*=0.036$, $t_*=T_9=0.044\dots=2/45$}
\label{zhshel:fig7}
\end{figure}

For parameters (\ref{zhshel:eq:3.10}) the numerical results are presented in figures~\ref{zhshel:fig6}, \ref{zhshel:fig7}.



\section{Interaction between strong and weak discontinuities}\label{zhshel:sec:7}

It should be clarified that the solution presented in section~\ref{zhshel:sec:6}  is valid before time $t=T_9$.
In moment $t=T_9$ the left boundary of the zone $\mathbb{Z}_9$ (line $x=x^1_w(t)$) intersects the right boundary of the zone  $\mathbb{Z}_2$ (line $x=x^1_s(t)$). It corresponds to the intersection of strong and weak discontinuities for the Riemann invariant $R^1$ (see figure~\ref{zhshel:fig4}).
Solving the equation $x^1_w(t)=x^1_s(t)$ one can obtain the interaction point $(X_9,T_9)$
\begin{eqnarray*}
T_9=\frac{\mu^2-q^1}{\mu^1-q^1}T_3, \quad X_9=x^1 + q^1 \mu^1 \mu^2 T_9.
\end{eqnarray*}

Starting at moment $t=T_9$ left boundary of zone $\mathbb{Z}_9$ is modified.
To get new boundary $x=\Phi(t)$ we use the  Rankine--Hugoniot conditions
(\ref{zhshel:eq:3.04}). The Riemann invariant $R^2$ as before is a constant ($R^2 = \mu^2$). It means that one of the Rankine-Hugoniot conditions (\ref{zhshel:eq:3.04}) is automatically satisfied, and another condition has the form
\begin{eqnarray}\label{zhshel:eq:7.02}
\frac{D}{\mu^2}\left[\frac{1}{R^1}\right]=
-\left[ R^1 \right].
\end{eqnarray}
On the new boundary $x=\Phi(t)$ of the zone $\mathbb{Z}_9$ (see figure~\ref{zhshel:fig4}) we have
\begin{eqnarray*}
R^1(\Phi(t)-0,t)=\mu^1, \quad R^1(\Phi(t)+0,t)=R^1_9(\Phi(t)+0,t),
\end{eqnarray*}
where $R^1_9(x,t)$ is the Riemann invariant $R^1$ for the zone $\mathbb{Z}_9$ that is obtained with the help of the relations  (\ref{zhshel:eq:6.07}), (\ref{zhshel:eq:6.08}).

We introduce the notation (compare with~(\ref{zhshel:eq:6.10}))
\begin{eqnarray*}
\rho^1(t)=R^1(\Phi(t)+0,t)=R^1_9(\Phi(t)+0,t).
\end{eqnarray*}

Using (\ref{zhshel:eq:6.04}), (\ref{zhshel:eq:6.08}), and (\ref{zhshel:eq:7.02}) we get
\begin{eqnarray}\label{zhshel:eq:7.05}
\frac{d\Phi(\beta)}{d\beta}=D=\mu^1\mu^2\rho^1(\beta),
\end{eqnarray}
\begin{eqnarray}\label{zhshel:eq:7.06}
\Phi(\beta)=\varphi(\tau) + \mu^2 \rho^1(\beta) \rho^1(\beta)(\beta-\tau),
\end{eqnarray}
\begin{eqnarray}\label{zhshel:eq:7.07}
\tau=t(\rho^1(\beta), \mu^2), \quad \varphi(\tau)=x(\rho^1(\beta), \mu^2).
\end{eqnarray}
Here $\beta$ is the current time $t$. We change the notation in order not to  confuse the current time $t$ and function $t(\rho^1,\mu^2)$.
We note that (\ref{zhshel:eq:7.06}), (\ref{zhshel:eq:7.07}) differ from (\ref{zhshel:eq:6.11})
only in that the parameter $\rho^1$ depends from time.
We also note that (\ref{zhshel:eq:7.05}) is not the equation for the function $\Phi(t)$ which is determined by (\ref{zhshel:eq:7.06}),  (\ref{zhshel:eq:7.07}). It is ODE for  the function  $\rho^1(t)$ which is written in the implicit form.

Adding the obvious initial conditions
\begin{eqnarray*}
\rho^1(T_9)=q^1, \quad \Phi(T_9)=X_9,
\end{eqnarray*}
\begin{eqnarray}\label{zhshel:eq:7.09}
T_9=\frac{\mu^2-q^1}{\mu^1-q^1}T_3, \quad X_9=x^1 + q^1 \mu^1 \mu^2 T_9,
\end{eqnarray}
we obtain the Cauchy problem (\ref{zhshel:eq:7.05})--(\ref{zhshel:eq:7.09}) for determination of the function $\rho^1(t)$.

Similar way one can study the interaction of discontinuities which occurs at the point $(X_{10},T_{10})$
\begin{eqnarray*}
T_{10}=\frac{\mu^1-q^2}{\mu^2-q^2}T_6, \quad X_{10}=x^2+q^2\mu^1\mu^2T_{10},
\end{eqnarray*}
when the lines $x=x^2_w(t)$ and $x=x^2_s(t)$ are intersected
(see figure~\ref{zhshel:fig4}). To obtain the solution it is enough to replace the indices $1 \rightleftharpoons 2$  and the functions  $\Phi \rightarrow \Theta$.

\begin{figure}[H]
\centering
\includegraphics[scale=0.49]{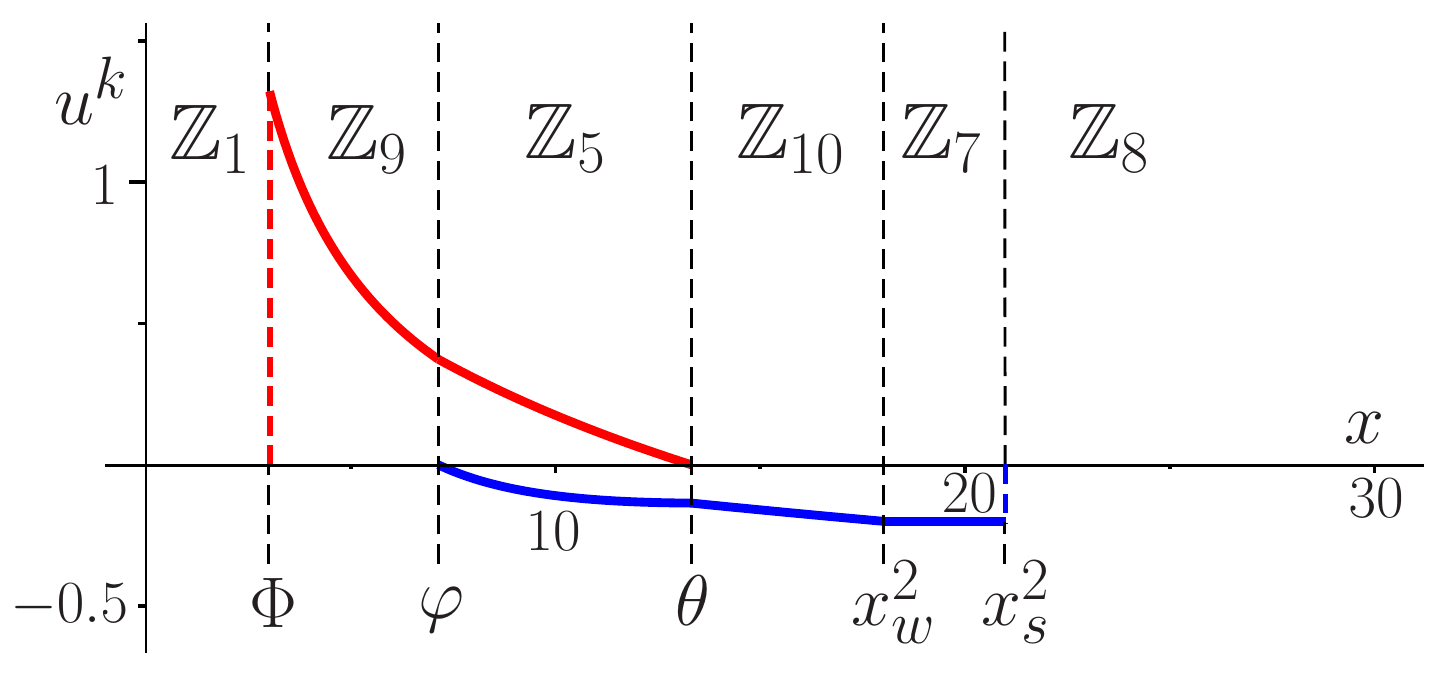}
\includegraphics[scale=0.49]{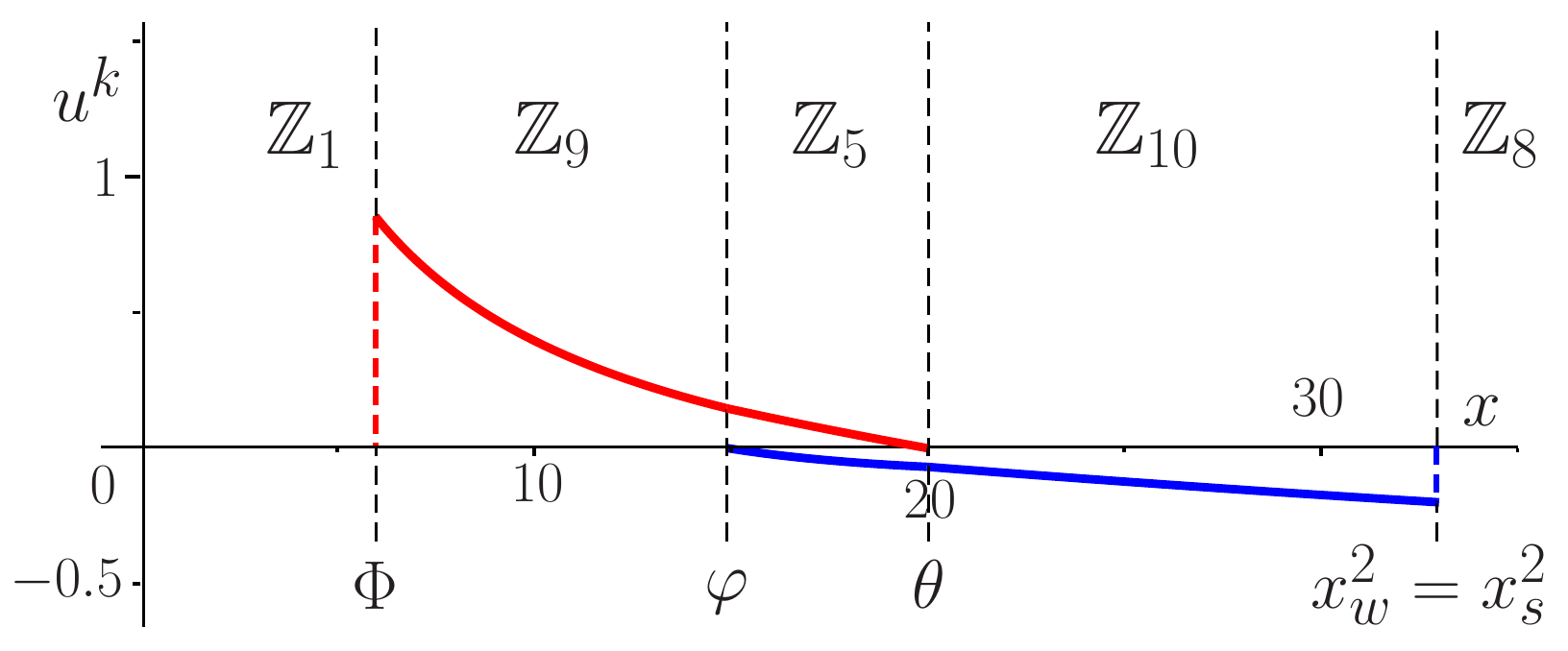}
\caption{Concentrations $u^k$ at the moments $t_*=0.05$, $t_*=T_{10}=0.08$}
\label{zhshel:fig8}
\end{figure}

\begin{figure}[H]
\centering
\includegraphics[scale=0.50]{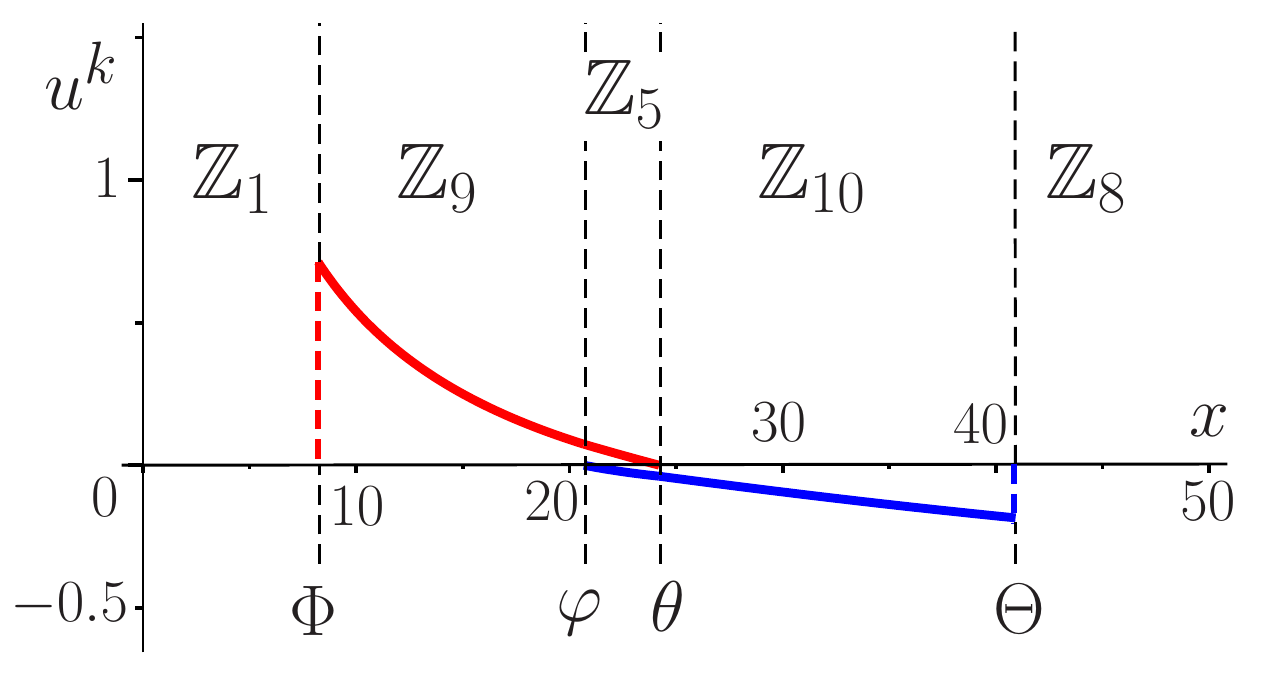}\qquad\qquad
\includegraphics[scale=0.50]{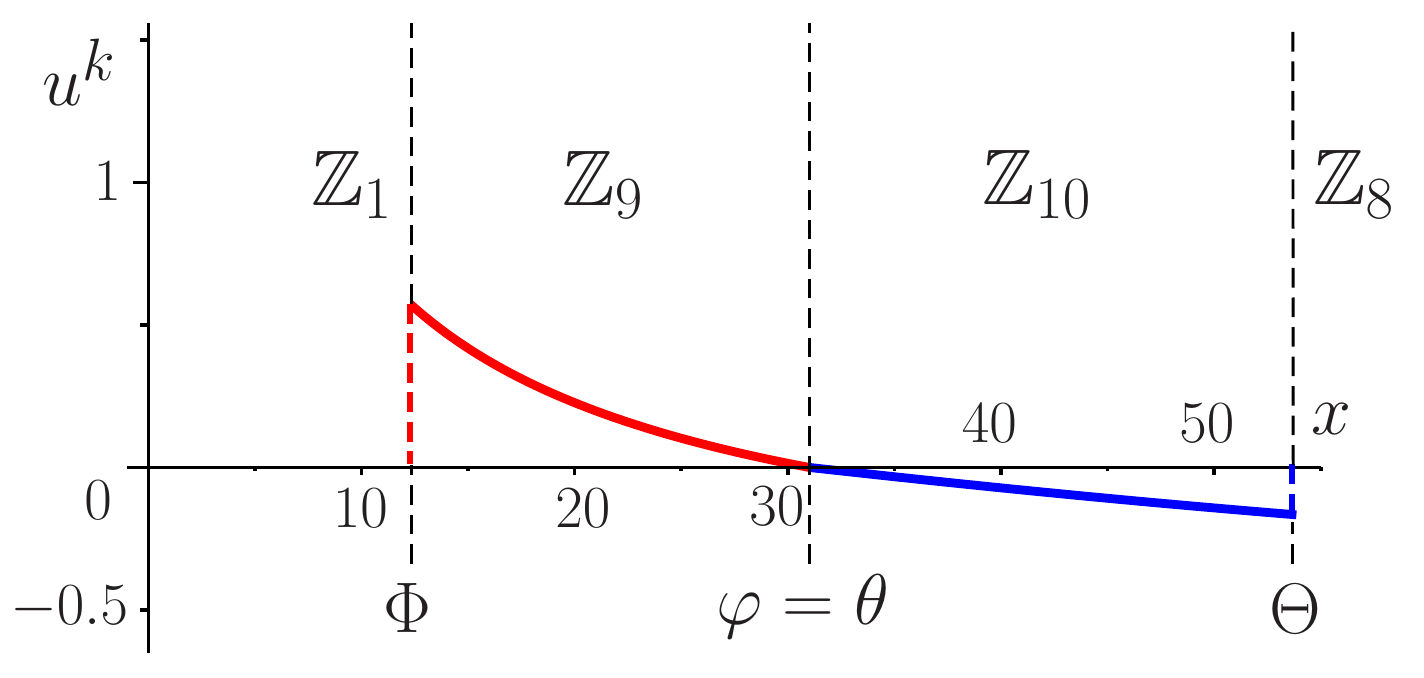}
\caption{Concentrations $u^k$ at the moments $t_*=0.1$, $t_*=T_{fin}=0.133\dots=2/15$}
\label{zhshel:fig9}
\end{figure}

In figures~\ref{zhshel:fig8}, \ref{zhshel:fig9}  the numerical results are presented for parameters (\ref{zhshel:eq:3.10}).




\section{The behavior of solutions at $t \geqslant T_{fin}$}\label{zhshel:sec:8}

The last interaction of the discontinuities is occurred in the point $(X_{fin},T_{fin})$ when the lines $x=\varphi(t)$ and $x=\theta(t)$
are intersected (see figure~\ref{zhshel:fig4}). Functions $\varphi(t)$, $\theta(t)$ are determined by
the relations (\ref{zhshel:eq:6.04}), (\ref{zhshel:eq:6.05}). The new zone $\mathbb{Z}_{11}$ appears (see figure~\ref{zhshel:fig4}):
\begin{eqnarray*}
\mathbb{Z}_{11}=\{(\mu^1, \mu^2), (x^1_f(t), x^2_f(t))\},
\end{eqnarray*}
where the boundaries of the zone are determined by 
\begin{eqnarray*}
x^1_f=X_{fin}+\mu^1 \mu^1 \mu^2 (t-T_{fin}), \quad x^2_f=X_{fin}+\mu^2 \mu^1 \mu^2 (t-T_{fin}).
\end{eqnarray*}

The boundaries of the zones $\mathbb{Z}_{9}$ and $\mathbb{Z}_{10}$ are also changed
\begin{eqnarray*}
\mathbb{Z}_9=\{(R^1(x,t), \mu^2), (\Phi(t), x^1_f(t))\}, \quad
\end{eqnarray*}
\begin{eqnarray*}
\mathbb{Z}_{10}=\{(\mu^1,R^2(x,t)), (x^2_f(t),\Theta(t))\},
\end{eqnarray*}
where functions $\Phi(t)$,  $R^1(x,t)$, $\Theta(t)$ and $R^2(x,t)$  can be obtained by the method which is described in section~\ref{zhshel:sec:7}.

The values $T_{fin}$, $X_{fin}$ are determined by the formulae (\ref{zhshel:eq:6.04}), (\ref{zhshel:eq:6.05}), where
$\rho^1=\mu^1$, $\rho^2=\mu^2$
\begin{eqnarray*}
T_{fin}=t(\mu^1,\mu^2), \quad X_{fin}=\varphi(T_{fin})=\theta(T_{fin})=x(\mu^1,\mu^2).
\end{eqnarray*}

We note that component concentrations in the zone $\mathbb{Z}_{11}$ are zero (see formulae (\ref{zhshel:eq:2.05}), where $R^1=\mu^1$ and $R^2=\mu^2$). Thus the mixture is separated into individual components. In other words, zone $\mathbb{Z}_{11}$  is  an `empty' zone  between the zones $\mathbb{Z}_{9}$ and $\mathbb{Z}_{10}$.

It is interesting that if  we a priori assume the separation of the mixture then
$T_{fin}$ can be obtained with the help of the formula (\ref{zhshel:eq:4.08}) without the results of sections~\ref{zhshel:sec:5}--\ref{zhshel:sec:7}
\begin{eqnarray}\label{zhshel:eq:8.05}
T_{fin}=T_{int}V(q^1,q^2|\mu^1,\mu^2).
\end{eqnarray}

\begin{figure}[H]
\centering
\includegraphics[scale=0.54]{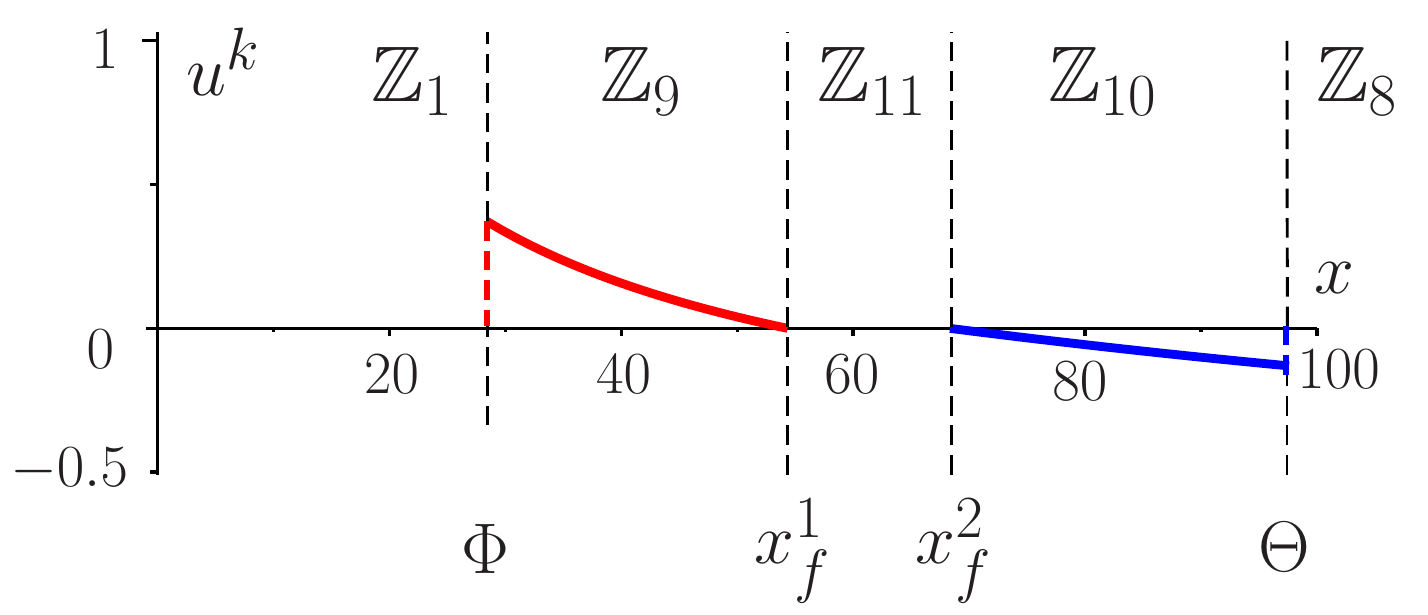}
\caption{Concentrations $u^k$ at the moment $t_*=0.25$}
\label{zhshel:fig10}
\end{figure}

For parameters (\ref{zhshel:eq:3.10}) the numerical results are presented in figure~\ref{zhshel:fig10}.


\section{Conclusion}

The existence of an analytic explicit formula for the Riemann function plays a key role for the method proposed.
The application of this method is not restricted to the problem of zone electrophoresis. Various equations of mathematical physics which can be studied by this method are listed for example in \cite{SenashovYakhno}. Such equations are the shallow water equations, gas dynamics equations, the Born-Infeld equation, nonlinear hyperbolic heat equation, Euler--Poisson--Darboux equation, chromatography equations etc. Note that there is an alternative way: this is application of the generalized hodograph \cite{Tsarev}.

If the initial data of the Cauchy problem allow to obtain an explicit formulae for the
commuting flows, then these flows can be used instead of the Riemann–Green function for
construction of the solution implicit form, i.\,e. $t(R^1,R^2)$, $x(R^1,R^2)$.

Another important result is the method  proposed for recovery of an explicit solution with the help of its implicit form.
In fact, it is the transformation of the Cauchy problem for PDEs into the Cauchy problem for ODEs on the level lines (isochrones).
This method allows not only to solve the Goursat problem, but also to solve effectively the Cauchy problem with arbitrary initial data (including discontinuous data, see \cite{Zhuk_Shir_ArXiv_2014_Part2}).

\begin{figure}[H]
\centering
\includegraphics[scale=0.6]{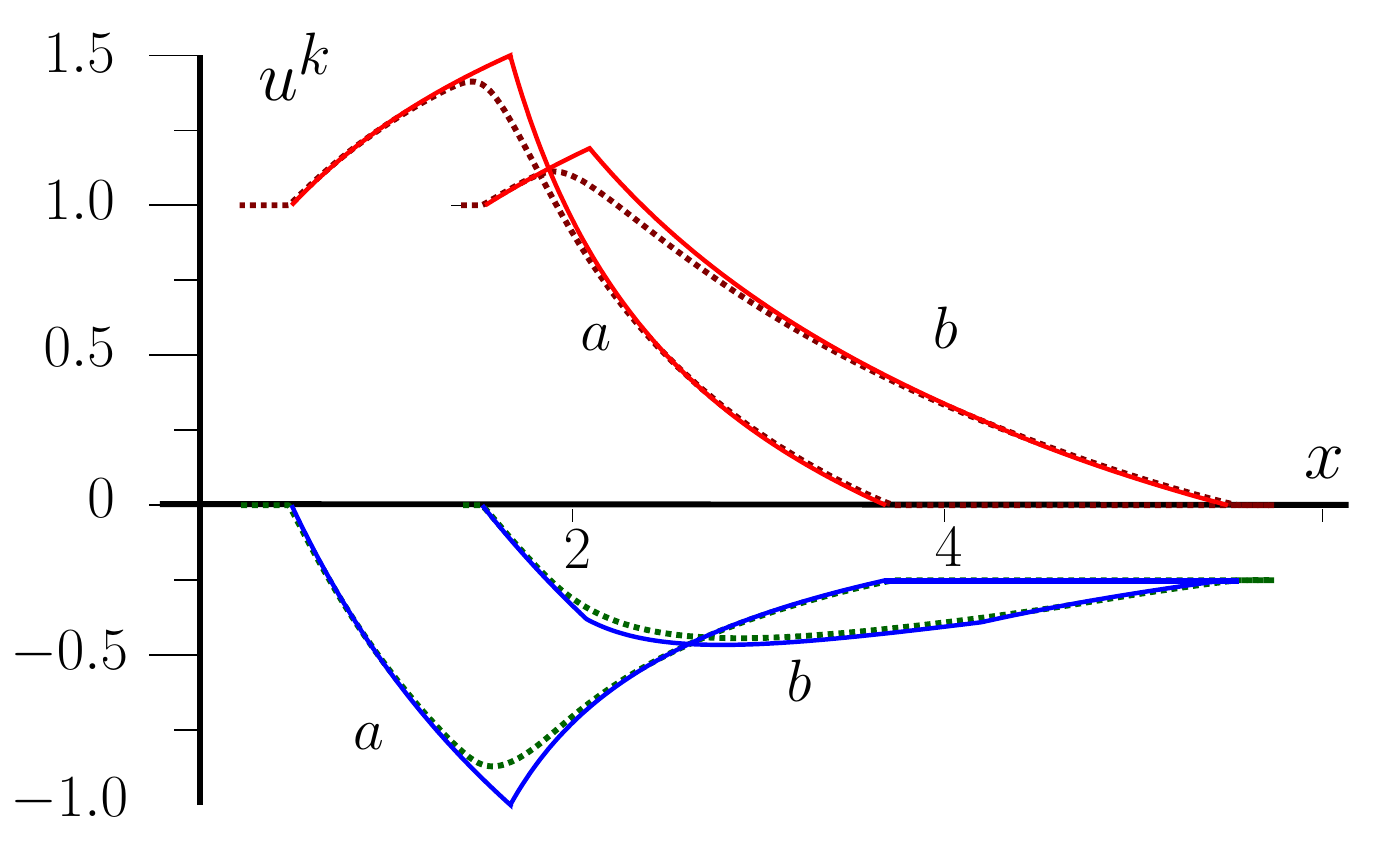}\\ [-2mm]    
\caption{Calculations (dotted lines)  and theory, $\left.a\right)$ $t=0.167$, $\left.b\right)$ $t=0.280$} 
\label{zhshel:fig11}
\end{figure}

Finally, to emphasize the importance of exact solutions, we present in figure~\ref{zhshel:fig11} the comparison of the exact solution with numerical solution.

The solid line corresponds to exact solution of the Goursat problem  obtained in sections \ref{zhshel:sec:5}, \ref{zhshel:sec:6}, and the dotted line corresponds to the solution constructed by a numerical method (FEM and VOF). Figure~\ref{zhshel:fig11}  illustrates the smoothing of the numerical solution. It indicates to uncertain identification of weak discontinuities and, as a result, the impossibility of  the discontinuities interaction checking.

\section*{}
\ack
The research was supported by the Base part of the project 213.01-11/2014-1 (subject 2016 year)
the Ministry of Education and Science of the Russian Federation, Southern Federal University.
The authors are grateful to N.M. Polyakova for proofreading the manuscript.

\appendix


\section{Numerical-analytical method for the Cauchy problem}\label{zhshel:sec:9}

We give a brief description of the numerical-analytical method, which allows to solve the Cauchy problem
(\ref{zhshel:eq:2.03}), (\ref{zhshel:eq:3.01}) with any initial data (in detail, see \cite{Zhuk_Shir_ArXiv_2014_Part2}). The implicit solution of the problem can be written in the following form (compare with (\ref{zhshel:eq:4.04}))
\begin{eqnarray}\label{zhshel:eq:A.01}
t= t(a,b)\equiv \frac{2(b-a)}{(r^1-r^2)^3}-\frac{(r^1+r^2)F(a,b)}{(r^1-r^2)^3}+ \frac{2 r^1 r^2 G(a,b)}{(r^1-r^2)^3},
\end{eqnarray}
where
\begin{eqnarray*}
  r^1=R^1_0(b), \quad r^2=R^2_0(a), \quad F(a,b)=
\int\limits_{a}^{b}f(\tau)\,d\tau,
\quad
G(a,b)=
\int\limits_{a}^{b}g(\tau)\,d\tau,
\end{eqnarray*}
\begin{eqnarray*}
f(\tau)=\frac{R^1_0(\tau)+R^2_0(\tau)}{R^1_0(\tau) R^2_0(\tau)}, \quad
g(\tau)=\frac{1}{R^1_0(\tau) R^2_0(\tau)}.
\end{eqnarray*}

On the $(a,b)$-plane we choose the point $(a_*,b_*)$ which identifies the isochrone $t_*=t(a_*, b_*)$.
In practice point $(a_*,b_*)$ can be selected with the help of the level lines of the function $t(a,b)$
for some interval of the parameters $a$, $b$.

We assume that the isochrone is parameterized by functions
$a=a(\mu)$, $b=b(\mu)$, where $\mu$ is the parameter. To obtain the functions $a=a(\mu)$, $b=b(\mu)$ on isochrone  we solve the Cauchy problem
\begin{eqnarray}\label{zhshel:eq:A.02}
\frac{da}{d\mu}=-t_b, \quad \frac{db}{d\mu}=t_a, \quad \frac{dF}{d\mu}=F_{\mu}, \quad \frac{dG}{d\mu}=G_{\mu},
\quad \frac{dX}{d\mu}=X_\mu,
\end{eqnarray}
\begin{eqnarray*}
a\bigr|_{\mu=0}=a_*, \quad b\bigr|_{\mu=0}=b_*, \quad
F\bigr|_{\mu=0}=F_*, \quad G\bigr|_{\mu=0}=G_*, \quad
X\bigr|_{\mu=0}=X_*,
\end{eqnarray*}
where the right parts of differential equations can be calculated in the explicit form
\begin{eqnarray*}
F_{\mu}=f(b)t_a+f(a)t_b, \quad G_{\mu}=g(b)t_a+g(a)t_b,
\end{eqnarray*}
\begin{eqnarray*}
X_\mu=(\lambda^2(r^1,r^2)-\lambda^1(r^1,r^2))t_a t_b, \quad X\equiv x(a,b).
\end{eqnarray*}
Integrating from $b=a_*$ to $b=b_*$ the Cauchy problem
\begin{eqnarray}\label{zhshel:eq:A.03}
\frac{dY(b)}{db}=(\lambda^2(r^1(b),r^2(a_*)) t_b(a_*,b), \quad Y(a_*)=a_*,
\end{eqnarray}
\begin{eqnarray*}
\frac{dF(a_*,b)}{db}=f(b), \quad \frac{dG(a_*,b)}{db}=g(b), \quad F(a_*,a_*)=0, \quad G(a_*,a_*)=0,
\end{eqnarray*}
we obtain values $X_*=Y(b_*)$, $F_*=F(a_*,b_*)$, $G_*=G(a_*,b_*)$.

For any value of the parameter $\mu$ the solution of the original problem (\ref{zhshel:eq:2.03}), (\ref{zhshel:eq:2.04}) with initial data (\ref{zhshel:eq:3.01}) has the form
\begin{eqnarray*}
R^1(x,t_*)=R^1_0(b(\mu)), \quad R^2(x,t_*)=R^2_0(a(\mu)), \quad x=X(\mu).
\end{eqnarray*}

The method described can be used for the solution of any Cauchy problem, if the explicit relation for the function $t(a,b)$ is given. Strictly speaking, `numerical part' of the method is required only to solve the Cauchy problems (\ref{zhshel:eq:A.02}) and (\ref{zhshel:eq:A.03}). In some cases, these problems can be solve explicitly, as, for instance, for the Born--Infeld equation (see e.g. \cite{Zhuk_Shir_ArXiv_2014_Part3}). In a sense, the method is a generalization of the method of characteristics for the case of two equations. In this connection see paper \cite{CurroFuscoManganaro}, which develops quite similar method for solving equations of chromatography.


\section*{References}
\begin{thebibliography}{99}
\small

\bibitem{RozhdestvenskiiYanenko}
Rozdestvenskii B L, Janenko N N 1983 {\it Systems of Quasilinear Equations and Their Applications
to Gas Dynamics (Translation of Mathematical Monographs vol 55)} (Provience, RI: American  Mathematical Society)

\bibitem{ElaevaIzvestiya}
Elaeva M S 2010 Interaction between strong and weak  discontinuities in the Riemann problem for hyperbolic equations
{\it University news. North-Caucasian Region. Natural Sciences series} {\bf 6} 14--19




\bibitem{Copson}
Copson E T 1958 On the Riemann-Green Function {\it Arch. Ration. Mech. Anal.} {\bf 1} 324--348


\bibitem{Bizadze}
Bitsadze A V 1981 {\it Some Classes of Equations in Partial Derivatives} (Moscow: Nauka) p~448



\bibitem{SenashovYakhno}
Senashov S I, Yakhno A 2012 Conservation laws, hodograph transformation and boundary value problems of plane plasticity
{\it SIGMA} {\bf 8} 16


\bibitem{FerapontovTsarev_MatModel}
Ferapontov E V, Tsarev S P 1991 Systems of hydrodynamic type that arise in gas chromatography. Riemann invariants and exact solutions
{\it Mat. Model} {\bf 3} 82--91



\bibitem{BabskiiZhukovYudovichRussian}
Babskii V G, Zhukov M Yu, Yudovich V I 1989 {\it Mathematical Theory of Electrophoresis} (New York: Plenum Publishers Corporation) p~241


\bibitem{Zhuk_Shir_ArXiv_2014_Part2}
Shiryaeva E V, Zhukov M Yu 2014 {\it Hodograph Method and Numerical Integration  of Two Quasilinear Hyperbolic Equations. Part~II. The Zone
Electrophoresis Equations} arXiv:1503.01762



\bibitem{ElaevaMM}
Elaeva M S 2010 Investigation of zone elecrophoresis for two component mixture
{\it Mat. Model} {\bf 22} 146--160

\bibitem{Elaeva_Diss}
Elaeva M S 2011 {\it Mathematical modeling of capillary zone electrophoresis. Dissertation of candidate of physico-mathematical sciences:
05.13.18} (Rostov-on-Don: SFU) p~161

\bibitem{Elaeva_ZhVM}
Elaeva M S 2012 Separation of two component mixture under action an electric field
{\it Comp. Math. and Mat. Phys} {\bf 52} 1143--1159

\bibitem{Lax}
Lax P D 2006 {\it Hyperbolic Partial Differention Equations} (Provience, RI: American  Mathematical Society)

\bibitem{Liu}
Liu T-P 2000 {\it Hyperbolic and Viscous Conservation Laws} (CMBS-NSF Regional Conference Series in
Applied Mathematics, SIAM)



\bibitem{Pavlov_Preprint}
Pavlov M V 1987 {\it Hamiltonian formalism of electophoresis equations. Integrable equations of hydrodynamics}
(Moscow: Preprint of Landau Inst. for Theor. Phys. 17) p~7



\bibitem{Tsarev}
Tsarev S P 1991 The geometry of hamiltonian systems of hydrodynamic type. The generalized hodograph method
{\it USSR Izvestiya Mathematics} {\bf 54:5} 1048--1067


\bibitem{Zhuk_Shir_ArXiv_2014_Part3}
Shiryaeva E V, Zhukov M Yu 2015 {\it Hodograph Method and Numerical Integration  of Two Quasilinear Hyperbolic Equations. Part~III. Two-Beam Reduction of the Dense Soliton Gas Equations} arXiv: arXiv:1512.06710




\bibitem{CurroFuscoManganaro}
Curro C, Fusco D, Manganaro N 2015 Exact description of simple wave interactions in multicomponent chromatography
{\it J. Phys. A: Math. Theor.} {\bf 48} 015201


\end{thebibliography}
\end{document}